\documentclass[12pt]{article}
\title{The noise in the circular law and the Gaussian free field}
\date{\today}
\author{Brian Rider
  and B\'alint
Vir\'ag
}

\usepackage{amsfonts}
\usepackage{graphicx}
\usepackage{amsmath}
\usepackage{amsthm}
\usepackage{amssymb}

\theoremstyle{theorem}
    \newtheorem{theorem}{Theorem}
    \newtheorem{lemma}[theorem]{Lemma}
    
    \newtheorem{corollary}[theorem]{Corollary}

\theoremstyle{definition} % For roman text in the body

\theoremstyle{remark} % For an italic header, more subtle than definition style

\newcommand\be{\begin{equation}}

\newcommand\ee{\end{equation}}

\newcommand{\PP}{\mathcal P}

\newcommand{\comment}[1]{}
\newcommand{\Z}{{\mathbb Z}}
\newcommand{\UU}{{\mathbb U}}
\newcommand{\R}{{\mathbb R}}
\newcommand{\CC}{{\mathbb C}}
\newcommand{\DD}{{\mathbb D}}

\newcommand{\ev}{\mbox{\bf E}}

\newcommand{\Var}{\mbox{\rm Var}}

\newcommand{\Cov}{\mbox{\rm Cov}}
\newcommand{\sm}{{\raise0.3ex\hbox{$\scriptstyle \setminus$}}}

\newcommand{\re}[1]{(\ref{#1})}

\newcommand{\al}{\alpha}

\newcommand{\su}{\twoheadrightarrow}
\newcommand{\ep}{\varepsilon}
\newcommand{\ra}{\rightarrow}

\newcommand{\ld}{\lambda}
\newcommand{\dl}{\delta}

\newcommand{\zn}{z^{\oplus n}}
\newcommand{\zm}{z^{\oplus m}}

\begin{document}
\maketitle
\begin{abstract}
Fill an $n \times n$ matrix with independent complex Gaussians of
variance $1/n$. As $n\to \infty$, the eigenvalues $\{ z_k\}$ converge to a sum
of an $H^1$-noise on the unit disk and an independent $H^{1/2}$-noise
on the unit circle. More precisely, for $C^1$
functions of suitable growth,
the distribution of $\sum_{k=1}^n (f(z_k)-\ev f(z_k))$
converges to that of a mean-zero Gaussian with variance given by
the sum of the squares
of the disk $H^1$ and the circle $H^{1/2}$ norms of $f$.
%$\frac{1}{4 \pi} \|f\|_{H^1}^2+  \frac{1}{2} \|f\|_{H^{1/2}}^2$
As a consequence, with $p_n$ the characteristic polynomial,
it is found that $\log
|p_n|-\ev \log|p_n|$ tends to the planar Gaussian free field
conditioned to be harmonic outside the unit disk. Further, for
polynomial test functions $f$,  we prove that the limiting
covariance structure is universal for a class of models including
Haar distributed unitary matrices.
\end{abstract}

\section{Introduction}

\begin{figure}
\centering
\includegraphics[height=1.5in]{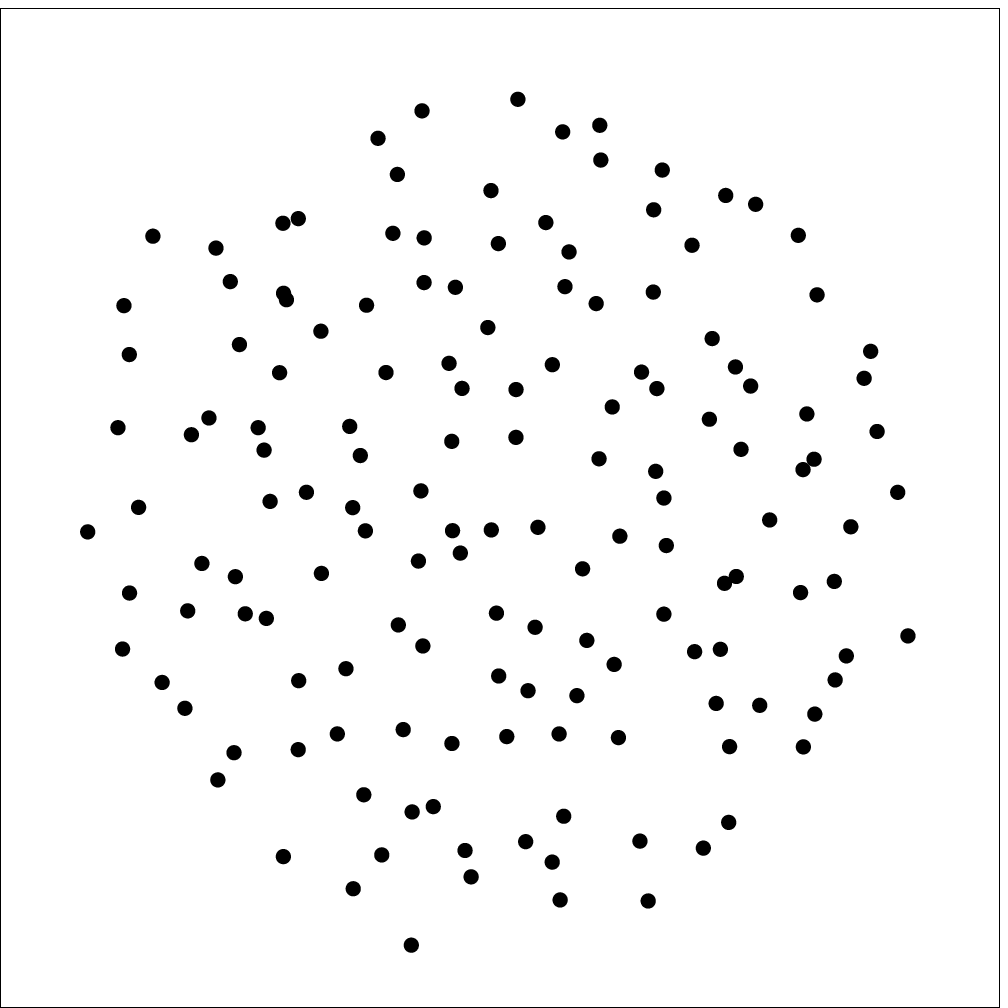}\hspace{1in}
\includegraphics[height=1.5in]{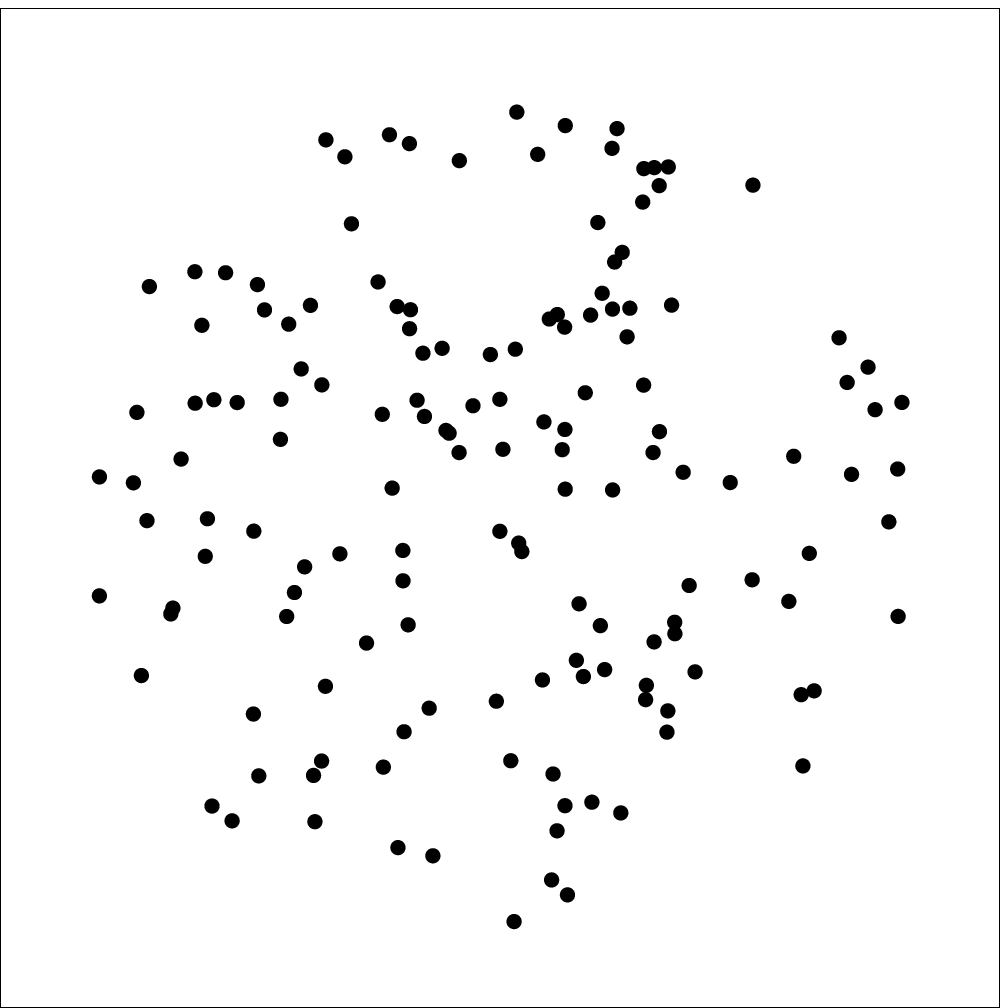}
\caption{\label{f1} Ginibre eigenvalues and uniform points in
$\UU$, $n=150$.}
\end{figure}

Consider the {\em Ginibre ensemble}, that is the $n
\times n$ random matrix in which all entries are independent complex
Gaussians of mean zero and variance $1/n$. The eigenvalues $z_1, z_2, \dots, z_n$
form a
point process, a realization of which is depicted on the left of Figure \ref{f1}. A first glance at
this picture suggests that the eigenvalues are uniformly
distributed in the unit disk $\UU$. Indeed, Bai \cite{Bai97} has shown that, with probability one,
the empirical measure converges to the uniform distribution on $\UU$.
This is the circular law of the title.

\begin{figure}
\centering
\includegraphics[height=1.5in]{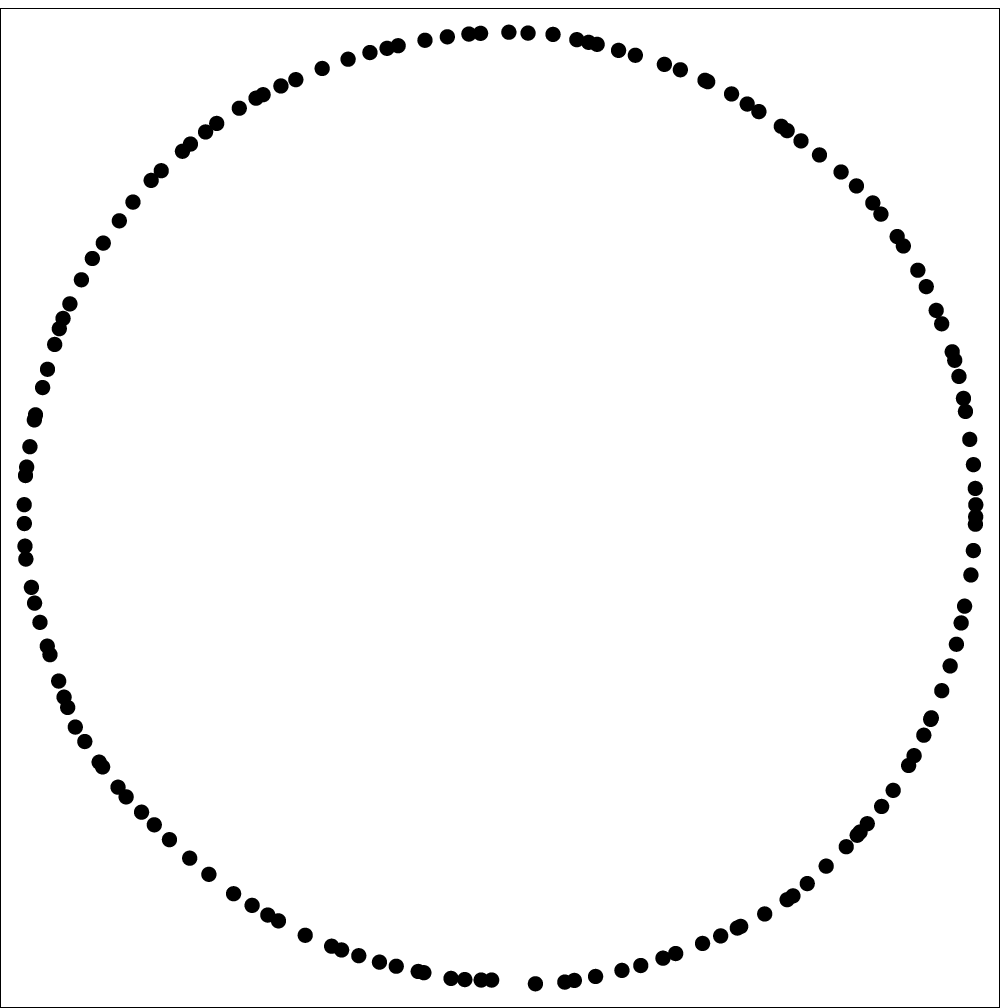}\hspace{1in}
\includegraphics[height=1.5in]{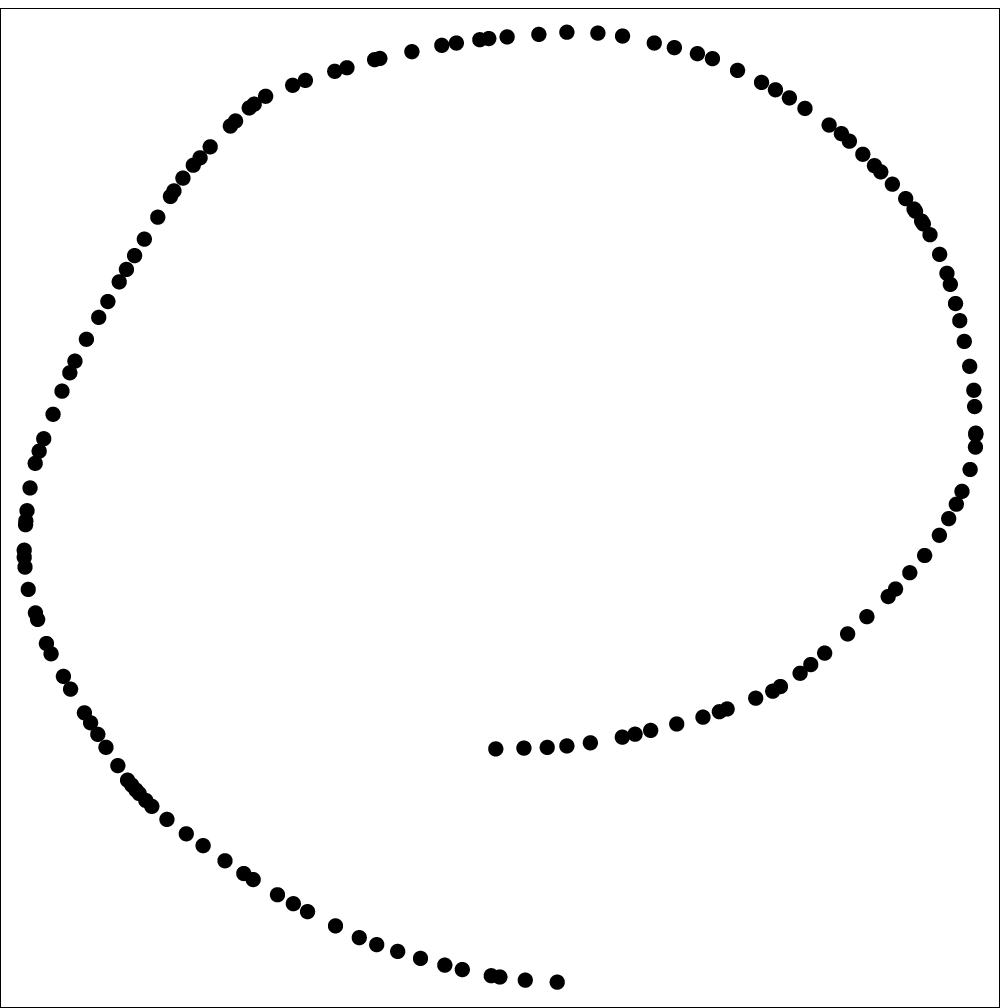}
\caption{\label{f1.5} The running eigenvalue sums $\sum_{k=1}^m
z_k$, $m =1, \dots, 150$,  in which the $z_k$ are ordered by their
arguments, and the same for the independent points.}
\end{figure}

Compared to $n$ points dropped independently and uniformly on $\UU$
(see the right of Figure \ref{f1}), the Ginibre points  are clearly
more regular. For example, their sum, the matrix trace, is complex
Gaussian with variance 1, while for the  independent points the
variance is $\frac{2}{3}n$ (see Figure \ref{f1.5} for comparison).
Our main theorem addresses this phenomenon in greater generality,
rigorously verifying a prediction of Forrester \cite{For99}.
\begin{theorem}[Noise  limit for eigenvalues]
\label{main} Let $f: \CC \ra \R$ possess continuous partial derivatives
in a neighborhood of $\UU$, and grow at most  exponentially
at infinity.  Then, as $n\to\infty$, the distribution of the
random variable $\sum_{k=1}^n \left(f(z_k) - \ev f(z_k)\right)$
converges to a normal with variance $\sigma^2_f+\tilde \sigma^2_f$,
where
$$\sigma^2_f=  \frac{1}{4\pi}\| f \|_{H^1(\UU)}^2 =\frac{1}{4\pi} \int_{\UU} | \nabla f |^2 d^2
z$$ is the squared Dirichlet ($H^1$) norm of $f$ on the unit disk,
and
$$\tilde \sigma^2_f=\frac{1}{2}\| f
\|_{H^{1/2}(\partial \UU)}^2 = \frac{1}{2}\sum_{k \in \Z}  |k|  |
\hat{f}(k) |^2$$ is the $H^{1/2}$-norm on the unit circle $\partial
\UU$, where $ \hat{f}(k) = \frac{1}{2 \pi} \int_{0}^{2 \pi}
f(\theta) e^{-i k \theta} d \theta$ is the $k$-th Fourier
coefficient of $f$ restricted to $|z| = 1$.
\end{theorem}

We mention that the centralizer $ \sum_{k=1}^n  \ev f(z_k)$ is easily replaced
by $ n \times \frac{1}{\pi} \int_{\UU} f(z) d^2 z$ with no change to the outcome.
More interesting is the covariance structure that appears above,
being composed of two independent terms: one for the boundary and one for
the bulk. An interpretation of this limiting  noise  is provided upon
considering the characteristic polynomial
$p_n(z) = \prod_{k=1}^n (z-z_k)$ which has fluctuations
described in terms of  the Gaussian free field.

The planar Gaussian free field (GFF)  is a model which  has received
considerable recent attention as the scaling limit of uniformly
random (discrete) $\R^2\to \R$ surfaces, though it apparently
 has  not previously been connected with any matrix model.
 While
the realizations of the GFF fluctuate too wildly to allow it to
be defined pointwise,  it may be defined as a random
distribution $h$ as follows: for all suitable test functions $f$,  the random
variable $\langle f, h \rangle _{H^1(\CC)}$ is Gaussian with variance
$\|f\|^2_{H^1(\CC)} = \int_{\CC} |\nabla f|^2 \, d^2 z$.  In more classical language this is the
Gaussian Hilbert space for $H^1(\CC)$, identifying the GFF as
 the most natural
2-dimensional analogue of Brownian motion.

As the Ginibre eigenvalues accumulate in the unit disk,  $\log |p_n(z)|$
tends to be harmonic away from the unit disk as $n \to \infty$.
Thus  the following version of the GFF appears naturally in the
limit.

\begin{corollary}[The Gaussian free field limit]
\label{gffcor} For the Ginibre characteristic polynomial $p_n(z) = \prod_{k=1}^n  (z - z_k)$, let
$$h_n(z)= \log |p_n(z)|-\ev \log |p_n(z)|.$$
Then $h_n$ converges weakly
without normalization to $h^*$, the
planar Gaussian free field conditioned to be harmonic outside the
disk. More precisely, for functions $f$ as above we have
$$
\int h_n(z)f(z) \,d^2z \Rightarrow \int h^*(z)f(z) \,d^2z
$$
in distribution, and the same holds for the joint distribution for
the integrals against multiple test functions $f_1,\ldots,f_k$.
\end{corollary}

Next, by the Cram\'er-Wold device, a version of Theorem
 \ref{main} holds for complex valued test functions $f$ of like regularity
and growth. In this case the $H^1$ noise term is replaced by the $L^2$
norm of the $\bar{\partial} = (1/2)( \partial_x + \partial_y)$ derivative
of $f$.  Then, if it is further assumed that $f$ is complex
analytic in a neighborhood of
$\overline \UU$, the covariance structure simplifies significantly:
\[    \sum_{k=1}^n \left(f(z_k) - \ev f(z_k)\right)
 \, \Rightarrow \,   \mathcal{ N}_{\CC} \Bigl(  0 ,  \frac{1}{\pi} \| f^{\prime} \|_{L^2(\UU)}^2
        \Bigr)=
\mathcal{ N}_{\CC} \Bigl(  0 , \|f\|_{H^{1/2}(\partial U)} \Bigr),
\]
where $\mathcal{ N}_{\CC}$ denotes complex Gaussian distribution
with given mean and variance. This should be compared with the
result of Rider and Silverstein \cite{RidSilv05}, where the same
central limit theorem is proved for random matrices with more
general than ${\mathcal N}_{\CC}(0,1)$ entries, but where the test
functions need to be analytic throughout $|z| \le 4$.

Perhaps surprisingly, exactly the same limiting covariance
structure arises from eigenvalues of Haar distributed unitary
matrices (Diaconis and Evans, \cite{DiacEv01}). While the
difference of the two limiting noises seems very clear (one lives
on the disk and one on the circle), analytic functions do not see
this. We show that this is not a coincidence by proving a
universality result.

The important point is that both the Ginibre eigenvalues and those
of Haar $U(n)$ comprise   {\it symmetric polynomial projection
(SPP) processes} in $\CC$. An SPP process is a determinantal
process defined by a pair $(\mu,n)$, where $\mu$ is a rotationally
invariant probability measure on $\CC$ and $n$ is the number of
points, as follows. For $k=n$ the joint intensity $p(z_1,\ldots,
z_k)$ with respect to $\mu^k$ is given by $\det(K(z_i,z_j)_{1\le
i,j\le k})$, for $K$  the integral projection kernel to the
subspace of polynomials of degree less than $n$ in $L^2(\CC,\mu)$.
In fact,  as soon as  this holds for $k=n$,  it holds for any $k$,
see for example \cite{BenKrPerVir05}.

\begin{theorem}[Universal noise limit for analytic polynomial test functions]
\label{secondthm} Consider a sequence of SPP processes defined by
$(\mu_n,n)$. Assume that for all integers $m$ we have
$$
M(n,2n+2m)/M(n,2n)\to 1,\qquad \mbox{where }M(n,k)= \int_{\CC}
|z|^{k}d\mu_n(z),
$$
as $n \to \infty$.   Then,
 for $f$ any polynomial in $z$,
\[
   \sum_{k=1}^n \left(f(z_k) - \ev f(z_k)\right)
        \, \Rightarrow\,   \mathcal{ N}_{\CC}
        \Bigl(  0 ,  \frac{1}{\pi} \| f^{\prime} \|_{L^2(\UU)}^2  \Bigr),
\]
in distribution.
\end{theorem}

The proof of Theorem \ref{secondthm} relies on the close connection
between SPP processes and Schur functions. For unitary matrices,
this connection is a consequence of Frobenius duality in
representation theory (see \cite{DiacEv01}); in the last part of the
paper we show that they are useful in a more general setting. For
the Ginibre ensemble $\mu_n=\mathcal{ N}_{\CC} (0,1/n)$, while for
unitary eigenvalues $\mu_n$ is uniform measure on the unit circle.
When $\mu_n$ is uniform measure in the unit disk we get the
(truncated) Bergman ensemble, for which the $n\to\infty$ limit is a
{\it discrete} point process in $\UU$ that agrees with the zeros of
the random power series with i.i.d.\ complex Gaussian coefficients,
see \cite{BenKrPerVir05}.  Each of these ensembles satisfies the
moment condition laid out in the statement, and in each case the
result actually holds for $f$ analytic.

\bigskip

{\bf Methods.} Central limit theorems of the above type have been
studied in the random matrix theory community for matrix ensembles
Hermitian ensembles, e.g., \cite{Jo98, Guion02, BaiSilv04,
AndZeit06}, and ensembles connected with the classical compact
groups, e.g., \cite{Jo97, Sosh00, DiacEv01, Wie02}. In all of these
cases, the eigenvalues form a one-dimensional point process either
on the real line or on the unit circle.

The situation for two-dimensional eigenvalue processes is more
complicated; indeed, it seems that this paper is the first to prove
a central limit theorem in two dimensions for general (not
necessarily analytic) test functions.
For the one-dimensional case, variants Wigner's method of moments
argument has been refined and extended to give central limit
theorems in several settings. However, this argument does not work
for non-analytic test functions, as such functions are not
expressible as a trace of a matrix polynomial.
While analyticity in one dimension is essentially a smoothness
assumption, in two-dimensions it is much more restrictive. This
makes the proof of the general circular law \cite{Bai97}  more
difficult than that of the semicircle law, and explains why the
corresponding general central limit theorem is still open (again,
see \cite{RidSilv05} for
analytic test functions).

Our proof is based on the determinantal structure of the point
process and a cumulant formula for determinantal processes
introduced by Costin-Lebowitz \cite{CosLeb95} and generalized by
Soshnikov \cite{Sosh00a,Sosh00,Sosh02}. In the present setting, it allows us
to establish a connection between joint cumulants of monomial test
functions in $z$ and $\bar z$ and a class of combinatorial objects
which we call rotary flows.

When the associated matrix integrals can be evaluated by separate methods, then
such a connection can be used to count intractable combinatorial
objects or evaluate difficult combinatorial sums arising, for
example, in statistical physics. Here the cumulants become
combinatorial sums over rotary flows. Fortunately,
we are able to establish (Section \ref{subcomb})
conditions under which these
combinatorial sums vanish as $n \ra \infty$, without explicitly computing
them for  any finite $n$ case. As it turn out,
this is sufficient to conclude the  asymptotic joint normality
of  monomial statistics.
This result is then lifted from polynomial to general $C^1$
functionals by a concentration estimate, yielding Theorem \ref{main}. The sum of these steps  occupies
Sections \ref{sec:cumulants} through \ref{sec:conc}. Theorem 3 is proved
in Section \ref{sec:analytic}, while, after some preliminaries in the next section, the connection
to the GFF is  detailed in Section \ref{sec:gff}.

The arguments introduced here work in the general rotational
invariant setting of SPP processes. However, beyond analytic test
functions, the structure of the limit depends intimately on the
process in question. In this paper we concentrate on the Ginibre
ensemble.

\section{Facts about the Ginibre ensemble}

%{\bf The Ginibre ensemble. }

The non-Hermitian matrix ensemble
of independent complex Gaussian entries
is named for Ginibre who
discovered \cite{Gin65} that the joint density of eigenvalues, $z_1$
through $z_n$, is given by
\begin{equation}
\label{dens}
   d\mu_n(z_1)\cdots d\mu_n(z_n) \;\frac{1}{{\cal Z}_n}\; \prod_{j < k} | z_j -
  z_k|^2.
\end{equation}
Here $d \mu_n = \frac{n}{\pi} e^{-n |z|^2} d^2 z$ is complex
Gaussian measure of variance $1/n$, and ${\cal Z}_n < \infty$ is a
normalizing constant.

A description equivalent to (\ref{dens}) is to say that the
eigenvalues of the Ginibre ensemble make up the  determinantal process
tied to the projection onto the subspace of  polynomials of
degree less than $n$ in $L^2(\CC,\mu_n)$;  the projection kernel with
respect to $\mu_n$ being
\[
   K_n(z, \bar w) =  \sum_{\ell = 0}^{n-1} \frac{n^{\ell}}{\ell !} z^{\ell} {\bar w}^{\ell}.
\]
That is, all finite dimensional correlation functions, or joint
intensities, are expressed in terms of determinants of the kernel
$K_n$: for $B_1, \dots, B_k$ any mutually disjoint family of Borel
subsets of $\CC$, \be \label{equ:det}
  \ev \Bigl[  \prod_{\ell=1}^k  \# \{ k:z_k  \in B_{\ell} \}   \Bigr]  =  \int_{  \prod_{\ell =1}^k B_{\ell}}
  \det \Bigl(  K_n(z_i, z_j)  \Bigr)_{1 \le i, j \le k} d \mu_n(z_1) \cdots d \mu_n(z_k).
\ee
One then refers to the process $(K_n, \mu_n)$.
Further background information on determinantal processes
may be found in \cite{BenKrPerVir05} and \cite{Sosh00a}. Immediate
from either (\ref{dens}) or (\ref{equ:det}) is that the
eigenvalues repel each other; Figure 1 provides a  somewhat
dramatic snapshot of this fact.

As the dimension tends to infinity, the reader is invited to check by a simple
calculation that
the mean density of Ginibre eigenvalues,
\[
 \ev \Bigl[ \frac{1}{n} \sum_{k=1}^n \delta_{z_k}(z) \Bigr] =  \frac{1}{n} K_n(z, \bar z) d \mu_n(z),
\]
tends  (weakly) to the uniform measure  on the the
unit disk $\UU$.   This is a particularly simple instance of the circular law
about which we describe  fluctuations.

\bigskip

\section{The Gaussian free field}\label{sec:gff}

For a domain $\DD \subset \CC$, let $H^1=H^1(\DD)$ denote the Hilbert
space completion of smooth functions with compact support in $\DD$
with respect to the gradient inner product $ \langle
f,g\rangle_{H^1(\DD)}   = \langle f,g\rangle_{H^1(\DD)} =
 \int_{\DD} \nabla f  \cdot \nabla g  \,  d^2 z$.
%(To be precise, the elements of the Hilbert space are equivalence
%classes containing functions that differ by a constant).
Note that the $H^1$-norm is invariant under conformal
transformations and complex conjugation.
%In particular, if $f\in H^1(\DD)$, then $g(z)=f(1/z)\in
%H^1$ of the transformed domain $\DD^{\prime}$,
%and $\|f\|_{H^1(\DD)}=\|g\|_{H^1(\DD^{\prime})}$. The same holds for
%$g(z)=f(1/\bar z)$, by invariance under reflection about the $x$ axis.

The Gaussian free field (GFF) on $\DD$ is defined as the Gaussian
Hilbert space connected with this norm. It may also be interpreted
as the random distribution $h$ such that for all $f\in H^1(\DD)$ the
random variable $\langle f,h\rangle_{H^1(\DD)} $ is centered normal
with
$$
\Var \left( \langle f,h\rangle_{H^1(\DD)} \right)= \| f
\|^2_{H^1(\DD)}.
$$
For more detailed information concerning the GFF, the paper
\cite{Sheff} is recommended.

In contrast, the $H^1$-noise on $\DD$ may be defined as the Gaussian
random distribution $\tilde h$ for which $\langle f,\tilde
h\rangle_{L^2}$ is centered normal and
$$
\Var\left(\langle f,\tilde h\rangle_{L^2} \right) = \| f\|^2
_{H^1(\DD)}.
$$
The two are related by noting that
$$\langle f,h\rangle_{H^1} =
%\langle\grad f,\grad h\rangle_{L_2} =
\langle - \triangle
f,h\rangle_{L^2}= \langle f, - \triangle h \rangle_{L^2}.
$$
That is, formally there is the identity in law $\tilde h =  -
\triangle h.$

\subsection*{Planar GFF  harmonic outside $\UU$ and
Corollary 2}

The space $H^1(\CC)$ decomposes into three orthogonal parts
(see \cite{Sheff}):  $ H^1(\CC) =  A_I  \oplus A_O \oplus A_H$.
 $A_I$ is the closure (in $H^1$) of smooth
functions supported in a compact set of the open unit disk.
$A_O$ is the similar closure of functions supported in a
compact set in the complement of the closed unit disk.  Last, $A_H$ is the
subspace spanned by continuous functions which are harmonic on the
complement of the unit circle.

We can now define the GFF $h$ on $\CC$ conditioned to be harmonic exterior
to $\UU$ as the projection $\PP_{IH}$ of the standard planar GFF to $A_I \oplus A_H$ in
the distributional sense: for $f\in H^1$,  set $\langle f,
\PP_{IH}h \rangle = \langle \PP_{IH}f, h\rangle$.  Note that for
$f\in H^1$ and continuous,  $\PP_{H} f$ is
the harmonic extension of the values of $f$ on the unit circle to
the whole plane. That is,  the resulting function is harmonic at infinity, or
$(\PP_{H}f)(1/z)$ is harmonic near the origin.

The proof of Corollary 2 requires one preliminary observation.

\begin{lemma}\label{normslemma} For  $\PP_{IH}$ the projection unto the  $A_I
\oplus A_H$ it holds that
$$
\|\PP_{IH}f \|_{H^1(\CC)}^2 = \|f\|_{H^1(\UU)}^2 +  \frac{1}{2}
\|\PP_H f\|_{H^1(\CC)}^2.
$$
For $f$ continuous in the neighborhood of $\partial U$ the second
term satisfies
\be
 \label{norm2}
 \frac{1}{2} \|\PP_H f\|_{H^1(\CC)}^2 =
 \|\PP_H f\|_{H^1(\UU)}^2 = \pi \| f \|_{H^{1/2}(\partial \UU)}^2.
 \ee
\end{lemma}

\begin{proof}
If $f\in H^1(\CC)$, there exists a function $g \in A_O$ so that
$s= f + g$ is symmetric with respect to inversion, {\em i.e.},
$s(z)=s(1/\bar z)$.  With $ \| \cdot \|$ denoting $H^1(\CC)$
unless specified otherwise, we  write
$$
\|s\|^2= \|\PP_{I}s\|^2+\|\PP_{O}s\|^2+\|\PP_{H}s\|^2.
$$
By symmetry and the conformal and conjugation invariance of $H^1$,
%(see \cite{Sheff})
we have that $ \|\PP_{I}s\| = \|\PP_{O}s\|, $ and also
$\|s\|^2_{H^1(\UU)}=\|s\|^2_{H^1(\CC\setminus \UU)}=  \|s\|^2/2$.
Hence,
$$
\|\PP_{I}s\|^2 = \|s\|_{H^1(\UU)}^2-  \frac{1}{2} \|\PP_Hs\|^2,
$$
%and since $g=s-f\in A_O$,
and we conclude that
$$
\|\PP_{IH}f\|^2 =\|\PP_{IH}s\|^2=\|s\|_{H^1(\UU)}^2 +
\frac{1}{2} \|\PP_Hs\|^2=\|f\|_{H^1(\UU)}^2 +  \frac{1}{2} \|\PP_H f\|^2  %\hfill \qed
$$
since $g=s-f\in A_O$. This proves the first claim.

The first equality in \re{norm2} follows from symmetry and
conformal (inversion) and conjugation invariance. Note that all
norms in \re{norm2} depend only on the values of $f$ on $\partial
U$, so we may assume that $f$ is harmonic in $\UU$, and drop the
projection. Harmonic functions in $\UU$ are spanned by the real
and imaginary parts of $z^k$, so it suffices to check that for two
elements $g_1,g_2$ of this spanning set
$$\langle g_1, g_2 \rangle _{H^1(\UU)} = \pi \langle g_1,g_2  \rangle_{H^{1/2}(\partial
\UU)},
$$ which is a simple computation.
\end{proof}

\begin{proof}[Proof of Corollary \ref{gffcor}]
%Of course,  $\|\PP_H f\|_{H^1(\CC)}^2 =  2 ||  \PP_{\UU} f
%||_{H^{1}(\UU)}^2$ where $\PP_{\UU} f$ denotes the harmonic
%extension only to the interior of the disk.
%%For $f\in H^1$ and continuous,  $\PP_{H} f$ is
%the harmonic extension of the values of $f$ on the unit circle to
%the whole plane (so that it is also harmonic at $\infty$, i.e.,
%$(\PP_{H}f)(1/z)$ is harmonic near the origin).
%Moreover, for $g$
%harmonic in $\UU$ and continuous in a neigborhood of $\partial \UU$,
%one easily checks that check that
%$\|g\|_{H^1(U)}^2 = {2\pi}\|g\|_{H^{1/2}(\partial \UU)}^2$.  For example,
%$g(z)$ for $|z| < 1$ may be expressed via the Poisson integral of its boundary data.
%The kernel is expanded in series, the gradient computed term-wise, and upon squaring
%and integrating over the unit disk, the series for ($2 \pi  \times$) the $H^{1/2}$-norm  squared  is recognized.
By Lemma \ref{normslemma} the limiting variance in Theorem
\ref{main} can be written as
$$ \frac{1}{4\pi} || f ||_{H^1(\UU)}^2 +
                  \frac{1}{2} || f ||_{H^{1/2}(\partial \UU)}^2
                  = \frac{1}{4\pi} \| \PP_{IH}f  \|
^2_{H_1(\CC)}.
$$
%In Theorem \ref{main}, what we have found for the Ginibre ensemble
%is an $H^1$ noise in the interior of the unit disk.
Now recall  $p_n(z) = \prod_{k=1}^n (z - z_k)$, the Ginibre
characteristic polynomial, and set
$$
\tilde h_n(z)=  \frac{1}{2 \pi} \log |p_n(z)| -  \frac{1}{2 \pi} \ev
\log |p_n(z)| .
$$
Since $\frac{1}{2 \pi} \triangle    \log |z| = $ the delta measure
at the origin in distribution, $ h_n(z) =  \triangle \tilde h_n(z)$
equals the centered counting measure of eigenvalues in the same
sense.  Thus,  by Theorem 1, we have that for any $f$ which is once
differentiable in  a neighborhood of $\UU$, the random variable
$\langle f,\tilde h_n \rangle_{L^2(\UU)} $ is asymptotically normal
with variance $\frac{1}{4\pi} \| \PP_{IH}f\| ^2_{H^1(\CC)}$ as $n
\ra \infty$. In this way we can identify the distributional limit of  $\tilde h_n(z)$
as the planar Gaussian free field,  conditioned to be harmonic outside
the unit disk.
\end{proof}

\section{Cumulants, polynomial statistics and rotary flows}
\label{sec:cumulants}

For any real-valued random variable $X$ the cumulants, ${\cal C}_k(X), k = 1, 2, \dots$,
are defined by the expansion
\be
\label{defcum}
     \log \ev [ e^{i t X} ] = \sum_{k=1}^{\infty}  \frac{(it)^k}{k!}  {\cal C}_k(X),
\ee
and carry information in the same manner as the  moments of $X$.  Important
here is the fact that the variable $X$ is
Gaussian
if and only if ${\mathcal C}_k(X) = 0$ for all $k > 2$.

Now, as it seems to have been first observed by Costin and
Lebowitz in \cite{CosLeb95}, the cumulants of linear statistics in
(any) determinantal point process have a particularly nice form.
Staying in $\CC$, choose a kernel $K$ and measure $\mu$ so that
$K$ defines a self-adjoint integral operator $\cal K$ on $L^2(\CC,
d \mu)$. Then, if  $\cal K$ is locally trace class with all
eigenvalues in $[0,1]$, $(K, \mu)$ determines a determinantal
point process (see \cite{BenKrPerVir05}, \cite{Mach75} or
\cite{Sosh00a}). That is, there is a process with all $k$-point
intensities satisfying the identity $(\ref{equ:det})$. With $X(g)
= \sum g(z_k)$ and $\{ z_k \}$ the points of $(K, \mu)$, the
$k$-th cumulant is written,
\begin{eqnarray}
\label{cum1}
{\mathcal C}_k (X(g)) & =  & \sum_{m=1}^k \frac{(-1)^{m-1}}{m}
     \sum_{k_1  +  \cdots + k_m = k \atop k_1 \ge 1,  \dots k_m \ge 1     }  \frac{k!}{k_1 !  \cdots k_m!}  \\
    &  &  \times \int_{\CC^k}    \Bigl( \prod_{\ell=1}^m  (g(z_\ell))^{k_\ell}  \Bigr)
                        K(z_1, \bar{z_2}) K(z_2, \bar{z_3}) \cdots K(z_m, \bar{z_1})
                         d \mu(z_1)  \cdots d \mu(z_m). \nonumber
\end{eqnarray}
Along with \cite{CosLeb95}, the above has been put to important use
in \cite{Sosh00}, \cite{Sosh02}.

For what we do here it will be convenient to cast the cumulants somewhat differently.
First let  $[k]=\{1,\ldots,k\}$.
For a function $\sigma : [k] \to [m]$, and $f \in
\mathcal G^k$, where $\mathcal G$ is an algebra
real-valued functions, define $\sigma f\in {\mathcal G}^m$ by
$$
(\sigma f)_j (z) = \prod_{i:\, \sigma(i)=j} f_i (z).
$$
Fix a functional
$\Phi_m : \mathcal G^m \to \R$ for each $1 \le m \le
k.$ For $f \in \mathcal G^k$ define
$$ \Upsilon_{k,m} (\Phi_m , f) = \sum_{\sigma[m]\su [k]} \Phi_m (\sigma f),
$$
where the sum is over all surjections (denoted by
$\su$). Further   define
\be
\label{Upsil}
\Upsilon_k (\Phi, f) = \sum^k_{m=1} \frac{(-1)^{m-1}}{m} \Upsilon_{k,m} (\Phi_m , f).
\ee
Then, for the  determinantal process $(K, \mu)$, we set
$$
\Phi_m (f_1 , \ldots, f_m ) =
  \int f_1 (z_1) \cdots f_m (z_m ) K(z_1 ,\bar  z_2 )K(z_2 ,  \bar z_3 ) \cdots K(z_m ,
\bar z_1)\, d \mu(z_1) \cdots d \mu(z_m),
$$
and replace the formula (\ref{cum1}) for the $k$-th cumulant
with the equivalent,
\be
\label{cum2}
\mathcal C_k (X(g)) =  \Upsilon _k (\Phi,(g, g, \ldots, g)).
\ee

This new expression, as a sum over surjections rather than over
ordered partitions, lends itself to the combinatorial approach we
take.  In the following, it is explained how the cumulant formulas,
(\ref{Upsil}) through (\ref{cum2}), simplify in the case that the
determinatal process in question has a rotation invariant reference
measure $\mu$ in $\CC$ and the linear statistic considered is a
(weighted) polynomial in $z, \bar z$.  The next section (Section
\ref{subcomb}) establishes when the limiting combinatorial sums
vanish. Finally, these facts are put together to compute the
asymptotics of the cumulants of polynomial statistics in the Ginibre
ensemble, establishing the CLT in that case.

\subsection*{Rotation invariance and rotary flows}

The following identities hinge on the rotation invariance of Ginibre
ensemble. To drive this point home, we set things up for general
determinantal point processes which share this feature.

Begin with a radially symmetric reference measure $\mu$ in the complex plane  normalized
so that $\mu(\CC) = 1$, and define the kernel,
\be
\label{Gker}
 K(z, \bar w) =
 \sum_{\ell =0}^L \lambda_{\ell} \, c_{\ell} \, (z \bar w)^{\ell},
 \ee
where it is assumed that  $\lambda_{\ell} \in [0, 1]$, $\int_{\CC}
|z|^{2\ell} d\mu = c_{\ell}^{-1} <  \infty$, while  $L = \infty$
is allowed.  Again, the results of \cite{Mach75} or \cite{Sosh00a}
(see also \cite{BenKrPerVir05}) imply that $(K, \mu)$ define
appropriate joint intensities.

Now return to the formula (\ref{cum2}) for the $k$-th cumulant
of the linear functional $X(g) = \sum g(z_k)$:
 $${\mathcal C}_k (X(g)) =
 \Upsilon _k (\Phi,(g, g, \ldots, g)).
 $$
 Since $\Upsilon(\Phi, (f_1 ,\ldots , f_k))$ is a
$k$-linear symmetric functional, it follows that if
$g$ is a polynomial,  $g = \sum a_k f_k$,
 then ${\mathcal C}_k (X(g))$ is a sum
of terms of the form $\Upsilon(\Phi, (f_1 ,\ldots ,
f_k))$ in the monomials $f_k$.
 Our goal is to understand the conditions which will make $\Upsilon_k
(\Phi, (f_1\,\ldots, f_k ))$ vanish (or vanish in the limit of some parameter,
soon to be the dimension for Ginibre).

Looking to take advantage of the radial symmetry,
we fix a simple test function $\zeta(z)$, a product of
$z$ and a radially symmetric weight function.
Then taking the monomials  $f_j$ to be  powers of
$\zeta$ and its conjugate,
$$
f_j(z) = (\zeta(z))^{\alpha_j}\overline{(\zeta(z))^{\beta_j}},
$$
we
expand the corresponding integral to find that
 $$
  \Phi_m (f_1 ,\ldots, f_m)=
  \sum_{q_1 ,\ldots,q_m=0}^n \,\int \prod_{j=1}^m
\zeta_j(z_j)^{\alpha_j}\overline{\zeta_j(z_j)^{\beta_j}}
\lambda j c _j (z_j \bar z_{j+1} )^{q_j} \,d \mu(z_1)
\cdots d \mu(z_m).
 $$

\begin{figure}
\centering
\includegraphics[height=2in]{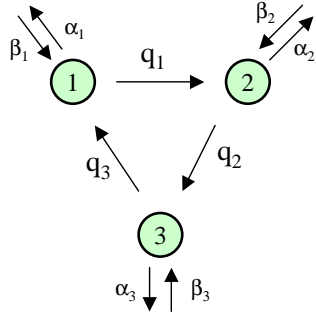}
\caption{\label{f2} A rotary flow: a cumulant term that does not
vanish}
\end{figure}

 Each term in the above sum can be thought of as a flow
network; the vertices are $1\ldots m$, and the directed edges are
$(j, j + 1)$, $j = 1 \ldots m$ (see Figure \ref{f2}). The exponent
of $z_j$ in the term is $\alpha_j + q_j$, and the exponent of $\bar
z_j$ is $\beta_j + q_{j-1}$; these have to agree or else the
integral over $z_j$ will vanish (granted by the rotation invariance
of the reference measure). In flow network terminology, the node law
has to hold (total inflow=total outflow), see figure \ref{f2}. This
implies that setting
\begin{eqnarray*}
 \gamma_j &=& \beta_j - \alpha_j,\qquad j=1,\ldots,m, \\
 \eta_j &=& \gamma_1 + \ldots + \gamma_j,
\end{eqnarray*}
and also $\ell=q_m$,
it is required that
\begin{eqnarray}
\label{conditions}
 q_j &=& \ell + \eta_j,\qquad \qquad  j=1,\ldots,m, \\
 \alpha_1 + \ldots + \alpha_m &=& \beta_1 +
 \ldots+\beta_m,  \nonumber
\end{eqnarray}
for the integral not to vanish.  We refer to $\gamma, \beta, \eta,q$
satisfying these conditions as a {\bf rotary flow}. Denoting
\begin{eqnarray*}
  M(q) &=&
  \int_\CC
  |z|^{q}  d \mu(z)
  = c_{q/2}^{-1},
 \\
  M(q,\kappa) &=&
  \int_\CC
  |z|^{q} |\zeta(z)|^{\kappa}
  d\mu(z) ,
%  \\
%  M(q,\kappa)&=&M(q,q,\kappa,\kappa),
\end{eqnarray*}
we find the following master formula,
 \begin{equation}
 \label{master}
 \Phi_m (f_1 , \ldots, f_m) =
  \sum^{n-1-\eta_{max}}_{\ell=-\eta_{min}}
  \prod^m_{j=1} \lambda_{\ell+\eta_j}
  \frac
   {M(2\ell + 2\eta_{j} -\gamma_j, \alpha_j+\beta_j)}
   {M(2\ell+2\eta_j)}.
 \end{equation}
 for the basic cumulant term.

In the flow language, the conditions (\ref{conditions}) say that
for each edge $(j,j+1)$ in the cycle we get a term
$$\lambda_{q_j}c_{q_j}=\lambda_{q_j}M(2q_j)^{-1}=\lambda_{\ell+\eta_j}M(2\ell+2\eta_j)^{-1}.$$
For the vertex $j$ the exponent of the $|z_j|$ factors
coming from the kernels $K(\cdot,\cdot)$ equals the
total flow within the cycle in or out of $j$:
$$
 q_{j-1}+q_j
 =2\ell + \eta_{j-1}+\eta_j
 = 2\ell + 2\eta_j -\gamma_j.
 $$ The total exponent of
$|\zeta(z_j)|$ is the total external flow in or out of
$j$, which equals $\alpha_j+\beta_j$.

%In the case $\zeta=z$ and $\lambda \equiv 1$,
%$q+\kappa=q'+\kappa'$ we have
%$$
%M(q,q',\kappa,\kappa') = a_q/a_{q+\kappa}.
%$$

\section{Combinatorial identities related to rotary flows}
\label{subcomb}

The main goal of this section is to establish general conditions
under which the combinatorial sums over rotary flows arising in
connection with the limits of joint cumulants vanish.

Let $V$ be a vector space; we will use $V=\R$ unless otherwise
specified.  For a function $\sigma : [k] \to [m]$, and $\alpha\in
V^k$ define $\sigma \alpha\in V^m$ by
\[
(\sigma\alpha)_j = \sum_{i: \, \sigma(i)=j} \alpha_i.
\]
Fix a function $\varphi_m : V^m \to \R$ for each $1\le m \le k$.
Define
\begin{equation*}
\Lambda_{k,m}(\varphi_m , \alpha) = \sum_{\sigma: [m]\su[k]}
\varphi_m(\sigma\alpha),
\end{equation*}
with the sum over all surjections, as well as
$$
\Lambda_k (\varphi, \alpha) = \sum^k_{m=1}\frac{(-1)^{m-1}}{m}
 \Lambda_{k,m} (\varphi_m , \alpha).
$$
The connection to the cumulant expressions above should be clear. By
a slight abuse of notation, we will use the shorthand $\varphi= x^p$
for the collection of functions $\varphi_m = \alpha_1^p + . . . +
\alpha^p_m$. Also let
$$
\Lambda_{k,m}(1)=\left|\{\sigma:[m]\su [k]\}\right|=\Lambda_{k,m}
(x, (1,0,\ldots,0)),$$
 and introduce the shorthand
$\Lambda_{k}(1)=\Lambda_{k} (x, (1,0,\ldots,0))$.

\begin{lemma}\label{linear}
For linear $\varphi$ it holds that,
$$
 \Lambda_k (x,\alpha) =
 \left\{
 \begin{array}{ll}
 \alpha_1, & \mbox{ if } k = 1, \\
 0,      & \mbox{ if } k \ge 2.
 \end{array}
 \right.
$$
\end{lemma}

\begin{proof}
The case $k=1$ is clear, so assume $k \ge 2$. Each $\varphi_k
(\alpha) = \alpha_1 + \ldots + \alpha_k$, which we denote by $s$. We
expand each $\Lambda_{k,m}$ according to whether $\sigma(k)$ has a
unique pre-image. This produces
$$
\Lambda_{k,m}(x, \alpha) = sm \Lambda_{k-1,m-1}(1) +
sm\Lambda_{k-1,m}(1),
$$
and so the expression for $\Lambda_k (x, \alpha)$ is just a
telescoping sum.
\end{proof}

\begin{lemma}
\label{quadratic} In the quadratic case we have
$$
 \Lambda_k (x^2,\alpha) =
 \left\{
 \begin{array}{ll}
 \alpha_1^2,  &  \mbox{ if } k = 1, \\
 2\alpha_1 \alpha_2,  & \mbox{ if } k = 2, \\
 0,      &  \mbox{ if }  k \ge 3.
 \end{array}
 \right.
$$
\end{lemma}

\begin{proof}
Again assume $k \ge 2$. Let
$\alpha^{-}=(\alpha_1,\ldots,\alpha_{k-1})$, and let $s^- = \alpha_1
+ \cdots + \alpha_{k-1}$. We expand $\Lambda_{k,m}$ in terms of
$\alpha_k$,  first considering the case where $\sigma(k)$ has a
unique pre-image $k$. Summing over values of $\sigma(k)$ we find
 $$
k\left(\Lambda_{k-1,m-1} (x^2 , \alpha^-) + \alpha_k^2
\Lambda_{k-1,m-1}(1)\right).
 $$
Now take the case when $\sigma(k)$ has more than one pre-image.
Expanding each quadratic term yields
$$
 m \Lambda_{k-1,m} (x^2 , \alpha^-) + m\alpha^2_k
 \Lambda_{k-1,m}(1) + \sum_{\sigma:[k]\su[m]}\sum_{i<k: \atop
 \sigma(i)=\sigma(k)} 2\alpha_k \alpha_i
$$
Performing the sum in the last term over the values of $\sigma(k)$
produces $2\alpha_k s^- \Lambda_{k-1,m}(1)$. It follows,
 \begin{eqnarray*}
 \Lambda_{k,m} (x^2 , \alpha) &=&
  m\Lambda_{k-1,m}(x^2 , \alpha^-) +
  m \Lambda_{k-1,m-1}(x^2 , \alpha^-) \\&+&
  m \alpha^2_k \Lambda_{k-1,m}(1 ) +
  m \alpha^2_k  \Lambda_{k-1,m-1} (1) \\ &+&
  2s^-\alpha_k \Lambda_{k-1,m} (1),
 \end{eqnarray*}
 and so,
 when we sum in $m$ to compute $\Lambda_k$, the first two rows
telescope leaving,
$$
 \Lambda_k (x^2 , \alpha) = 2s^-\alpha_k
 \Lambda_{k-1} (1).
$$
This equals $2\alpha_1 \alpha_ 2$ for $k = 2$, and $0$ for $k \ge 3$
by Lemma \ref{linear}.
\end{proof}

\begin{lemma}\label{combined}
Assume that each $\varphi_m$, $1\le m \le k$ is a (not necessarily
symmetric) quadratic polynomial in the $\alpha_i$. Let $b^{(i)}_m$
$i=0,1,2$ be the sum of the coefficients of the $\alpha_j^i$ terms
in $\varphi_m$; let $b^{(1,1)}_{m}$ denote the sum of the
coefficients of $\alpha_j\alpha_{j'}$ for $j\not=j'$. If
$$
b^{\lambda}_m=
 \left\{
 \begin{array}{ll}
 b^{(0)}, & \mbox{if } \lambda=(0), \\
 mb^{(i)}, & \mbox{if } \lambda=(1) \mbox{ or } \lambda=(2),\\
 m(m-1)b^{(1,1)}, & \mbox{if }  \lambda=(1,1), \\
 \end{array}
 \right.
$$
 for some $b^\lambda$ and all $m$, then
$$
 \Lambda_k (\varphi,\alpha) =
 \left\{
 \begin{array}{ll}
 b^{(0)}+b^{(1)}\alpha_1 + b^{(2)} \alpha_1^2, & \mbox{ if } k = 1, \\
 2(b^{(2)}-b^{(1,1)}) \alpha_1 \alpha_2, &  \mbox{ if } k = 2, \\
 0,      & \mbox{ if } k \ge 3.
 \end{array}
 \right.
$$
\end{lemma}

\begin{proof}
Denote by
$$
\tilde \varphi_m(\alpha)=\frac{1}{n!}\sum_{\sigma: [n]\su[n]}
\varphi_m(\sigma \alpha),
$$
the symmetrized version of $\varphi_m$. By definition
$\Lambda_{k,m}(\varphi_m,\alpha)$ is symmetric in the $\alpha_i$.
So, for $\sigma:[k]\su [m]$ and $y=\sigma \alpha$, if  we set
$s_\sigma=y_1+\ldots +y_m$, then
$s_\sigma=s=\alpha_1+\ldots+\alpha_k$ does not depend on $\sigma$ or
$m$.  Thus, $\sum_{i\not=j}y_iy_j=s^2-\sum_{i=1}^m y_i^2$, and
further
$$
\tilde \varphi_m(y) = b^{(0)}+b^{(1,1)}s^2+b^{(1)}\sum_{i=1}^m y_i+
(b^{(2)}-b^{(1,1)})\sum_{i=1}^m y_i^2.
$$
From here we also have
\begin{eqnarray*}
\Lambda_{k,m}(\varphi_m,\alpha)&=&\Lambda_{k,m}(\tilde
\varphi_m,\alpha) \\ &=&
 (b^{(0)}+b^{(1,1)}s^2)\Lambda_{k,m}(1)
 +b^{(1)}\Lambda_{k,m}(x,\alpha)+(b^{(2)}-b^{(1,1)})
 \Lambda_{k,m}(x^2,\alpha),
\end{eqnarray*}
and the claim follows from Lemmas \ref{linear} and \ref{quadratic}.
\end{proof}

\begin{lemma}\label{alphabeta}
Let now $V=\R^2$, and assume that each $\varphi_m: V^m \to \R$, is a
(not necessarily symmetric) quadratic polynomial in
$(a_i,\beta_i)_{i=1}^m$. Let $b^+_m$, $b^-_m$ denote the sum of the
coefficients of the $\alpha_i\beta_i$ and $\alpha_i\beta_j$,
$i\not=j$ terms in $\varphi_m$, respectively. Assume as well that
all other coefficients vanish, and that there exists $b^+$ and $b^-$
so that
$$
b^{+}_m=mb^+,
% \mbox{   and   }
\qquad b^{-}_m=m(m-1)b^-
$$
for all $m$. Then it holds
$$
 \Lambda_k (\varphi,\alpha) =
 \left\{
 \begin{array}{ll}
b^+ \alpha_1\beta_1, & \mbox{ if } k = 1, \\
(b^+-b^{-}) (\alpha_1\beta_2+ \alpha_2\beta_1), & \mbox { if } k = 2, \\
 0      & \mbox{ if }, k \ge 3.
 \end{array}
 \right.
$$
\end{lemma}

\begin{proof} By Lemma
\ref{combined} we may assume that each $\varphi_m$ is symmetric.
Setting $s_\beta=\beta_1+\ldots +\beta_k$, $s_\alpha=\alpha_1+\ldots
+ \alpha_k$, and for $\sigma:[m]\su[k]$ setting $y=\sigma \alpha$,
$w=\sigma\beta$, the arguments of Lemma \ref{combined} give
$$
\varphi_m ((y_1,w_1),\ldots,(y_m,w_m))= b^-s_\alpha s_\beta +
(b^+-b^-)\sum_{i=1}^m y_iw_i.
$$
We use polarization ($4y_iw_i=(y_i+w_i)^2-(y_i-w_i)^2$) for the
second term and find that,
$$
\Lambda_{k,m}(\phi,(\alpha,\beta)) = b^-s_\alpha
s_\beta\Lambda_{k,m}(1)+\frac{b^+-b^-}{4}\left(\Lambda_{k,m}(x^2,\alpha+\beta)-\Lambda_{k,m}(x^2,\alpha-\beta)
\right).
$$
The proof is finished by invoking Lemmas \ref{linear} and
\ref{quadratic}.
\end{proof}

\begin{lemma}[Spohn, Soshnikov] \label{etamax}
  Let $\varphi_m (\gamma_1 , \ldots, \gamma_m) =
\max(\eta_1 , \ldots , \eta_m)$, where $\eta_j = \gamma_1 + \ldots +
\gamma_{j}$. If $\alpha_1+\ldots + \alpha_k = 0$,  then
 $$
 \Lambda_k (\varphi, \alpha) =
 \left\{
 \begin{array}{ll}
 -|\alpha_1| = - |\alpha_2|, &  \mbox{ if } k = 2, \\
 0, &  \mbox{ if } k \neq 2. \\
 \end{array}
 \right.
$$
\end{lemma}

A result of this type was first discussed by Spohn in \cite{Spohn}
(though see \cite{Spitz} as well) which focusses on the
determinantal process on the line with the sine kernel
$\frac{\sin(\pi(x-y))}{\pi(x-y)}$. It also appears as the `Main
Combinatorial Lemma' in the the paper of Soshnikov, \cite{Sosh00},
where  it is used to track fluctuations of the eigenvalues of  the
Unitary group $U(n)$.
%However, the lemma follows easily from a
%paper by Kac, which relates $\varphi_m$ to the functional
%$$\varphi'_m(\gamma_1
%, \ldots, \gamma_m)=\sum_{k=1}^m \frac{\eta_k^+}{k} .$$ He shows
%that $\varphi'_m$ equals $\varphi_m$ after averaging over all
%permutation of the indices. In our notation
%$$\Lambda_{m,m}(\varphi_m,\eta)=\Lambda_{m,m}(\varphi'_m,\eta)$$

In our result, Lemma \ref{etamax} figures into the boundary
component of the limiting noise, and so may be thought of as the
root of the connection to  $U(n)$.

We finish with the note that in the formulation of \cite{Sosh00},
the $\eta'_j = \gamma'_1 + \ldots + \gamma'_{j-1}$.  That is, the
sum only goes up to index $j-1$).  This is  equivalent to the above
upon setting $\gamma'_i=-\gamma_{k+1-i}$,  so that
$\eta'_j=\gamma'_{1}+\ldots+\gamma'_{j-1}= -(\gamma_{k}+\ldots
+\gamma_{k+1-(j-1)})=\eta_{k+1-j}$ and  $\max(\eta'_1 , \ldots ,
\eta'_m)=\max(\eta_1 , \ldots , \eta_m)$.

\section{Polynomial statistics for the Ginibre ensemble}

We now apply the results of the previous two sections to the Ginibre
ensemble. That is, in (\ref{Gker}) we put
 $L = n < \infty$,
$\mu(z) = \mu_n(z) = \frac{n}{\pi} e^{-n |z|^2}$,
$c_k = \frac{n^k}{k!}$, $\lambda_k = 1$, and then take $n \ra \infty$.
Further, while the possibility of using weighted polynomials
may be important in general, in the case of Ginibre
it suffices to work with the naked monomials
$\zeta(z)=z$.

With these specifications,
$$
M(2\ell + 2\eta_j + \alpha_j-\beta_j, \alpha_j+\beta_j)= M(2\ell +
2\eta_j +2\alpha_j),
$$
and plainly
$$
\frac{M(2q+2\kappa)}{M(2\kappa)}=\frac{n^{q}(q+\kappa)!}{q!
n^{q+\kappa}}=n^{-\kappa}(q+1)\cdots (q+\kappa).
$$
Substituting into (\ref{master}), we have
\begin{equation}
\label{GinPhi}
\Phi_m( z_1^{\alpha_1}  {\bar z_1}^{\beta_1}, \dots, z_m^{\alpha_m}  {\bar z_m}^{\beta_m} )
=
\frac{1}{ n^s}\sum^{n-1-\eta_{max}}_
{\ell=-\eta_{min}} \prod_{j=1}^m (\ell + \eta_j +1)
\cdots{ (\ell + \eta_j + \alpha_j)},
\end{equation}
which for large $n$ reads,
$$
\frac{1}{ n^s}\sum^{n-1-\eta_{max}}_
{\ell=-\eta_{min}} \left[\ell^s+\ell^{s-1}\sum_{j=1}^m
\left[ \eta_j \alpha_j+{\alpha_j+1 \choose 2}\right]
+O(\ell^{s-2})\right],
$$
where again $s = \sum \alpha_i = \sum \beta_i$.
Next, use the fact that
$$
\sum_{\ell\,=\mbox{\,\scriptsize const.}}^n
\ell^s=\frac{n^{s+1}}{s+1}+\frac{n^s}{2}+O(n^{s-1}),
$$
valid for any $ s \ge 1$,  to put the
the asymptotics of (\ref{GinPhi}) into the form:
\begin{equation}
\label{GinPhi1}
\Phi_m(f_1, \dots, f_m)  =
\frac{n}{s+1}-(1+\eta_{\max}) + \frac{1}{2} +
\frac{1}{s}\sum_{j=1}^m  \left[ \eta_j \alpha_j+{\alpha_j+1
\choose 2}\right]+O(n^{-1}).
\end{equation}
Here  $ f_{\ell} = f_{\ell}(z, \bar z)$ is shorthand for $z^{\alpha_{\ell}} {\bar z}^{\beta_{\ell}}$.

For fixed $s$, the $O(1)$ term in (\ref{GinPhi1}) is a function of the
$\alpha_i$ and $\beta_i$, which we  denote
$\varphi_m$.  Summing the above into the basic cumulant term, we
find that
\begin{equation}
\label{GinPhi2}
\Upsilon_k(\Phi, (f_1 ,\ldots, f_k)) =
\frac{n}{s+1}\Lambda_k(1) +
\Lambda_k(\varphi,(\alpha,\beta))  + O(n^{-1}),
\end{equation}
and are in position to use the combinatorial tools developed
in Section \ref{subcomb}.

First of all, the coefficient of $n$ in (\ref{GinPhi2})
vanishes for $k\ge 2$ by Lemma $\ref{linear}$.
For the constant order term,
examining $\varphi_m$
shows it to be  a linear
combination of $\eta_{max}$, a quadratic polynomial in
$\alpha_i$ and a polynomial which is itself a linear
combination of the $\alpha_i\beta_j$.
Further, the sums of the
coefficients in the polynomial part
are as follows:
$$
\begin{array}{c|c|c|c|c|c}
 1&
 \alpha_i &
 \alpha_i^2 &
 \alpha_i\alpha_j, i\not=j &
 \alpha_i\beta_i &
 \alpha_i\beta_j, i\not=j \\
 \hline
 -1/2&
 m/(2s) &
 -m/(2s) &
 -m(m-1)/(2s) &
 m/s &
 m(m-1)/(2s)
\end{array}
$$
Now
combined,
Lemmas \ref{combined}, \ref{alphabeta} and
\ref{etamax} imply  that $\Lambda_k(\varphi, (\alpha, \beta))$,
and so the full cumulant,  vanishes in the limit
for $k > 2$.  In summary we have:

\begin{theorem}  For any real-valued polynomial $P(z, \bar z) = \sum c_{\alpha_k, \beta_k}
z^{\alpha_k} {\bar z}^{\beta_k}$, the Ginibre statistic $X_n(P) - \ev[ X_n(P)]$
converges in distribution to a mean-zero Gaussian as $n \ra \infty$.
\end{theorem}

As for the limiting variance, choose a pair of monomials
$z^{\alpha_1} {\bar z}^{\beta_1}$ and $z^{\alpha_2} {\bar
z}^{\beta_2}$ with $\alpha_1+\alpha_2=\beta_1+\beta_2$, and note
that our cumulant asymptotics enter the picture as in
 \begin{eqnarray}
 \label{secumu}
 \lim_{n \ra \infty}  \Cov (z^{\alpha_1} {\bar z}^{\beta_1} ,  z^{\beta_2} {\bar z}^{\alpha_2} )
 & = &   \lim_{n \ra \infty}  \Bigl(
  \Phi_1( z^{\alpha_1 + \alpha_2} \bar{z}^{\beta_1 + \beta_2} ) -
 \Phi_2( z^{\alpha_1} {\bar z}^{\beta_1} , z^{\alpha_2} {\bar z}^{\beta_2} )  \Bigr) \\
 &  = & \max(0, \beta_1 - \alpha_1)  +   \frac{\alpha_1 \beta_2}{s}.  \nonumber
  %\\
 %&&
 %\mathcal C(z^{\alpha_1}\bar z^{\beta_1},
 %            z^{\alpha_2}\bar z^{\beta_2})
 %  =\frac{|\,\alpha_1-\beta_1|}{2}
 %  +\frac{\alpha_1\beta_2+\alpha_2\beta_1}
 %  {\alpha_1+\alpha_2+\beta_1+\beta_2}.
 \end{eqnarray}
On the left, $z^{\alpha_1} {\bar z}^{\beta_1}$  and  $z^{\beta_2} {\bar z}^{\alpha_2} $
stand in for the full linear statistics, and  $\alpha_2,\beta_2$ switch rolls from left
to right  since covariance is conjugate linear in the second argument.
The limit itself is read off
from (\ref{GinPhi1});   why the corresponding covariance
vanishes in the limit if $\alpha_1+ \alpha_2 \neq \beta_1 + \beta_2$
should also be clear.

Next observe that
\be
\label{covt1}
   \frac{1}{ \pi }  \int_{\UU}  {\bar \partial}  (z^{\alpha_1} {\bar z}^{\beta_1})
        \overline{{\bar \partial}  (z^{\beta_2} {\bar z}^{\alpha_2})} d^2 z      =
         \frac{\alpha_2 \beta_1}{\alpha_1+\alpha_2}  =  \frac{\alpha_2 \beta_1}{s},
\ee
and
\be
\label{covt2}
   \sum_{k >0}  k  \, (z^{\alpha_1} {\bar z}^{\beta_1})^{\wedge}(k)   \, \overline{ (z^{\beta_2} {\bar z}^{\alpha_2})^{\wedge}(k)}  =
     \max(0,  \alpha_1 - \beta_1).
 \ee
A little algebra will show that the sum of (\ref{covt1}) and (\ref{covt2}) equals the final expression
in (\ref{secumu}).
By linearity we may conclude:
with any real-valued polynomials $f$ and $g$ in $z$ and $\bar z$,
$$
\lim_{n \ra \infty}
\Cov(f,g)=
\frac{1}{4\pi}\langle f,g \rangle_{H^1(\UU)} +
\frac{1}{2}\langle
f,g\rangle_{H^{1/2}(\partial \mathbb U)}.
$$
The general covariance/variance asymptotics occupy the next section.

\section{Concentration}
\label{sec:conc}

This section is devoted to the following estimate.

\begin{theorem}
\label{Conc} For linear statistics $X_n(f)$  and $X_n(g)$ in the Ginibre ensemble,
\be
\label{limitvar}
  \lim_{n \ra \infty}    {\Cov}( X_{n}(f), X_n(g)  )  =
      \frac{1}{\pi}  \int_{\UU}   {\bar \partial} f(z)     \overline{\bar \partial g}(z) \,   d^2 z   +
                \sum_{k > 0}  k   \widehat{f} (k)  \, \overline{\widehat{g}(k)},
\ee
as long as
$f$  and $g$ possess continuous partial derivatives
in  a neighborhood of $\UU$ and are otherwise bounded as in $ |f(z)| \vee |g(z)| \le C e^{c |z|}$.
\end{theorem}

Granting (\ref{limitvar}), the central limit theorem extends immediately from polynomials
of the form
\[
    P_m(  z, \bar{{z}}  )    =  \sum_{\alpha_k + \beta_k \le M}  c_{\alpha_k, \beta_k}
     z^{\alpha_k}  {\bar{z}}^{\beta_k},
\]
to general test functions $f(z)$ satisfying the growth and regularity conditions
assumed in the statement.

First, by the Stone-Weirstrauss theorem, we can find a polynomial
$P^{\prime}(z, \bar z )$ such that $\|  P^{\prime}(z, \bar z ) -
{\bar d} f \|_{L^{\infty}(|z| \le 1+ \delta)} \le \ep$ for
whatever $\ep >0$.  Certainly the $L^{\infty}$-norm (of the
derivatives) on the larger disk controls both the $H^{1}(\UU)$ and
$H^{1/2}( \partial \UU)$ norms, and there exists a sequence of
polynomials $P_m(z, \bar z)$ such that
\begin{equation}
\label{approx}
       \Bigl| \Bigl| f - P_m \Bigr| \Bigr|_{H^1(\UU) \cap H^{1/2}(\UU) } \ra 0,
\end{equation}
as the degree $m \ra \infty$.  At each step, $P_m$ is just the anti-derivative
(in $\bar z$) of the polynomial approximating ${\bar \partial} f$.

Next, denote the centered statistic by
\[
   Z_{n}(f) = X_{n}(f) -  \frac{n}{\pi} \int_{\UU} f(z) d^2 z.
\]
With $f$ as in the above result, a separate (and easy) calculation shows that
\[
     \ev [ X_{n}(f) ] -  \frac{n}{\pi} \int_{\UU} f(z) d^2 z = o(1),
\]
and it follows that the family $\{ Z_n(f) \}$ is tight as $n \ra
\infty$. Denote a given subsequential limit by $Z_{\infty}(f)$.
For polynomial test functions, we have proved that $Z_{n}(P_m)$
converges in distribution to a mean zero complex Gaussian,
$Z_{\infty}(P_m)$.  Now pick a smooth bounded $\phi: \CC \ra [0,
\infty)$ with $ \|\nabla \phi\|_{\infty} \le 1$,  and note that,
\begin{eqnarray*}
  \Bigl|  \ev  \phi ( Z_{\infty }(f) )  - \ev  \phi (Z_{\infty}(P_m) )  \Bigr|  & \le &   \limsup_{N \ra \infty}
   \ev  \Bigl|     \phi( Z_{n}(f) )  -  \phi(Z_{n}(P_m))           \Bigr|  \\
   & \le   &   \limsup_{n \ra \infty}  \ev [  | Z_{n}(f - P_m) |^2 ]^{1/2},
\end{eqnarray*}
where things have been fixed so the right hand side tends to zero as $m \ra \infty$.  This appraisal
is independent of the special subsequence chosen in the definition of $Z_{\infty}(f)$. Hence,
by appealing once more to Theorem \ref{Conc} and (\ref{approx}) the large $m$ limit of
$Z_{\infty}(P_m)$ identifies the limit
distribution of $Z_{n}(f)$ unambiguously and completes the proof of Theorem \ref{main}.

\subsection*{Verification of the limiting variance}

The opening move
in the proof of Theorem \ref{Conc} is
the covariance formula,
\begin{eqnarray*}
{\Cov} \bigr( X_n(f) , X_n(g)  \Bigr)  & =  & \int_{\CC}
                            f(z)   \, \overline
                            {g(z)}    \, K_n(z,\bar{z}) \, d \mu_n(z)  \\
                            &   & -
                             \int_{\CC}  \int_{\CC}
                             {f(z)} \, \overline{g(w) }  \,  |  K_n(z,\bar{w})  |^2   d \mu_n(z)  \, d \mu_n(w),
\end{eqnarray*}
valid for any determinantal point process.  Then, each appearance of
the test functions $f$ and $g$ is expanded by way of the dbar
representation. That is the fact that, for any $f$ once continuously
differentiable in a domain $\DD \subset \CC$ and continuous up to
the boundary $\partial \DD$, it holds,
\begin{equation}
\label{Dbar}
f(\zeta) =  - \frac{1}{\pi i}  \int_{\DD}  \frac{ {\bar \partial}  f(w) }{w -\zeta} \, d^2 w \
                + \ \frac{1}{2 \pi i} \int_{\partial  \DD}
                       \frac{f(w)}{w - \zeta}  \, d w,
\end{equation}
see for example \cite{Bell92}.

The formula (\ref{Dbar}) allows one to pass from the general
variance asymptotics to one for the particular functional $z \mapsto
(z - \cdot)^{-1}$.  The decay properties of the Ginibre ensemble
allow further simplifications. In particular, with $f(z)$ once
differentiable in  say $|z| < 1+\ep$, decompose as in  $f(z) = f(z)
\psi(z) + f(z) (1 -\psi(z))$ for a smooth $\psi$, equal to $1$ on
$\UU$ and vanishing for $|z| > 1+ \ep/2$. Next, for any $\ep > 0$
and $c > 0$,
\begin{eqnarray}
\label{ginbound}
     \int_{|z| \ge 1 + \ep } e^{c |z|^2 } K_n(z, \bar{z})  d \mu_N( z )
     & =  &
     \sum_{\ell=0}^{n-1}  \frac{n^{\ell+1}}{ \pi \ell !}  \int_{|z| \ge 1 + \ep }   e^{c |z|^2}  |z|^{2 \ell}  e^{-n |z|^2} d^2 z  \\
     & \le & n \times \frac{n^n}{(n-1)!} \int_{(1+\ep)^2}^{\infty} r^{n-1} e^{-(n-c) r}  \, dr  \le C^{\prime}  n e^{- n \ep^2/4},
                      \nonumber
\end{eqnarray}
by a simple Laplace estimate. From the growth assumption on $f$ and
noting $|K(z,\bar{w}) |^2 \le K(z,\bar{z}) K(w, \bar{w})$, it
follows that the variance of $X_n(f)$ agrees with that of $X_n{(f
\psi)}$ up to $o(1)$ errors. The same pertains to the covariance,
and it suffices to take the test functions compactly supported from
the start, ignoring the boundary integral in the dbar representation
(\ref{Dbar}). In short, the covariance may be read off from the $n
\ra \infty$ behavior of
\begin{eqnarray}
\label{Ginvar}
\lefteqn{  \frac{1}{\pi^2}
                 \int_{|\nu| < {1+ \ep}} \int_{|\eta| < {1+\ep}}
                        {\bar \partial}  f(\nu) \overline{ {\bar \partial} g (\eta)}  \,  \times } \\
 &  &  \left\{  \int_{\CC}  \frac{1}{z-\nu}   \, \overline
                            { \frac{1}{{z} - {\eta}}  }    \, K_n(z,\bar{z}) \, d \mu( z) -
                             \int_{\CC} \int_{\CC}
                                                          \frac{1}{z-\nu} \, \overline
                             {  \frac{1}{{w}- {\eta}} }   | K_n(z,\bar{w}) |^2  \,  d \mu_n( z)   d \mu_n( w)
                              \right\}     d^2 \nu    d^2 {\eta},  \nonumber
\end{eqnarray}
in which
$f(z)$ and $g(z)$ and their first partial derivatives are assumed to vanish along $|z| = 1+ \ep$.

Defining,
\begin{eqnarray}
\label{phidef}
   \phi_n(\nu,\eta)  & =  & \int_{\CC}
                            \frac{1}{z-\nu}   \, \overline
                            { \frac{1}{{z}- {\eta}}  }    \, K_n(z,\bar{z}) \, d \mu_n( z) \\
                            &         & -
                             \int_{\CC} \int_{\CC}
                             \frac{1}{z-\nu} \, \overline
                             {  \frac{1}{{w}- {\eta}} }  \, K_n(z,\bar{w})   K_n(w,\bar{z})  \, d \mu_n(z)  \, d \mu_n( w), \nonumber
\end{eqnarray}
the proof of Theorem \ref{Conc} splits into two parts.   For the interior (or $H^1$) contribution we have:

\begin{lemma}
\label{delta}
For fixed $\nu$ with $|\nu|< 1$,
\[
    \phi_n(\nu, {\bar \eta})  \ra \pi \, \delta_{\nu}(\eta),
\]
in the sense of measures as $n \ra \infty$.  It follows that,
\[
   \lim_{N \ra \infty}  \frac{1}{\pi^2}  \int_{|\nu| < 1}  \int_{| \eta| < 1+ \ep}  {\bar \partial f} (\nu)  \overline{ {\bar \partial f} (\eta)}  \phi_N( \nu, {\bar \eta} )
   \, d^2 \nu  d^2 \eta      =  \frac{1}{\pi}  \int_{\UU}  |{\bar \partial} f (\nu)|^2   d^2  \nu,
\]
for bounded and continuous  ${\bar \partial} f$.
\end{lemma}

For the boundary, (or $H^{1/2}$) contribution, we prove separately:

\begin{lemma}
\label{boundary}
With ${H_{+}^{1/2}(\partial \UU)}$ corresponding to the usual sum extended
only over positive indices, it holds that
\[
    \lim_{n \ra \infty}    \, \frac{1}{\pi^2} \int_{1 < |\nu|, |\eta| < 1+\ep}  {\bar \partial} f( \nu )
     \overline{ {\bar \partial} f (\eta ) } \, \phi_n(\nu, {\bar \eta} )
    \, d^2\nu \,   d^2 \eta  =
     \| f \|_{H_{+}^{1/2}(\partial \UU)}^2,
\]
whenever $f$ is continuously differentiable and vanishes on $|z| = 1 + \ep$.
\end{lemma}

Note that we have set $f=g$ in the above two statements.  This has
only been done for the sake of slightly more compact expressions.
No generality is lost: covariances can always be recovered from
variances.

\begin{proof} [Proof of Lemma \ref{delta}]
We begin by showing that,
\begin{equation}
\label{spike}
     \lim_{n \ra \infty}   \int_{\mathbb B}  \phi_n(\nu, {\bar \eta}) \, d^2 \eta   =  \pi,
\end{equation}
for $|\nu| < 1$ and $\mathbb B$ any disk about the origin containing $\nu$.
A residue calculation gives,
$$
    \int_{|\eta| < b}  \frac{ d^2 \eta}{z-\eta}  =
\left\{ \begin{array}{ll} \pi \bar{z} & \mbox{ for }  |z| < b \\
       \frac{\pi b^2}{z} & \mbox{ for }  |z| > b  \end{array} \right. ,
$$
and so, denoting the radius of ${\mathbb B}$ by $b> | \nu| = a$,
\begin{eqnarray}
\label{GS0}
   \int_{\mathbb B}   \phi_n(\nu, \bar \eta) \, d^2 \eta
      & =  & \pi   \int_{a <  |z| < b}   K_n(z, \bar{z}) \, d \mu_n(z) +
                         \pi  b^2 \int_{|z| > b}   \frac{1}{|z|^2}  K_n(z, \bar{z}) \, d \mu_n(z)  \\
       &     &    -   \pi  \sum_{\ell=0}^{n-2}   c_{\ell}  c_{\ell+1}
                         \int_{|z| > a}   {|z|}^{2\ell}    \, d \mu_n( z) \,
                         \int_{|w| <  b}  |w|^{2 \ell+2}  \, d \mu_n(w)   \nonumber \\
      &  &  - \pi b^2   \sum_{\ell=0}^{n-2}   c_{\ell}  c_{\ell+1}
                         \int_{|z| > a}   {|z|}^{2\ell}    \, d \mu_n( z) \,  \int_{|w| >  b}  |w|^{2 \ell}  \, d \mu_n(w), \nonumber
 \end{eqnarray}
after  expanding $(z-\nu)^{-1}$ in series where $|z|  > |\nu|$.
Recall that
$c_{k} = $  $c_k^n = \frac{n^{k}}{ k!}$.
Now recombine (\ref{GS0}) in
the form
\begin{eqnarray}
\label{GS1}
  \int_{\mathbb B}  \phi_n(\nu, \eta) d^2 \eta
  & = &   \pi  c_{N-1}  \int_{|z| > a} |z|^{2(n-1)}  \, d \mu_n(z) -  \pi
     \int_{|z| > b}  \Bigl( 1 -\frac{b^2}{|z|^2}  \Bigr) \, d \mu_n(z) \\
&   & - \,  \pi  \sum_{\ell = 0}^{n-2}   c_{\ell} c_{\ell+1}
                         \int_{|z| < a}   \int_{|w| > b}
                     |z|^{2\ell} |w|^{2\ell+2}   \Bigl( 1 - \frac{b^2}{|w|^2} \Bigr)    \, d \mu_n( z)
                        \, d \mu_n(w).             \nonumber
\end{eqnarray}
A simple computation shows that the first two terms of (\ref{GS1}) tend to $\pi$ and zero
respectively as $n \ra \infty$.  The final two terms of (\ref{GS1}) may in turn be bounded
by a constant multiple of
\begin{equation}
\label{probbound}
 \sum_{\ell = 0}^{n-2}   c_{\ell} c_{\ell+1}
                         \int_{|z| <a}   \int_{|w| > b}     |z|^{2\ell} |w|^{2\ell+2}    \, d \mu_n( z)
                        \, d \mu_n(w)
                         \le  \sum_{\ell = 0}^{n-2}  P \Bigl(   S_{\ell}^n \le a^2 \Bigr)
                                              P \Bigl( S_{\ell +1}^n \ge b^2   \Bigr)
\end{equation}
where $S_{\ell}^n$ denotes a sum of $(\ell+1)$ independent mean $1/n$ exponential random variables.
Standard large deviation estimates explain why
$P(S_{\ell}^n < a^2)$ is exponentially small in $n$ if $\ell > n( a^2 + \delta)$, with
an estimate of the same type holding for $P(S_{\ell}^n > b^2)$ and $\ell < n (b^2 -\delta)$.
Since there is a gap, {\em i.e.} $a < b$, each term in the above sum is of order $e^{-\gamma n}$
for  a $\gamma$ positive with $(b-a)^2$.  Since there are only $n$ terms, the verification of
(\ref{spike}) is complete.

It remains to show that we have decay away from $\eta = \nu$.  Given
$\delta > 0$,   we claim that
$$
         \lim_{n \ra \infty}    \int_{ {\mathbb A } \cap |\nu - \eta| > \delta }   \Bigl|  \phi_n(\nu, {\bar \eta})  \Bigr|  \, d^2 \eta  =  0,
$$
with any  $\mathbb A$ supported in $|z| \le 1 + \ep$.

First  fix $\eta$ as well as $ \nu $,  and take   $a = |\nu| \le |\eta| = b$ without
any loss.
Two types of expressions emerge from performing the integrations over $z$ and $w$ in the definition
of $\phi_n(\nu, {\bar{\eta}})$, an inner term corresponding
to $|z| < a$ and $|w| < b$, and an outer term corresponding to $|z|  > a $ and $|w|>  b$:
\[
    \phi_n(\nu, {\bar{\eta}}) =   {\phi}_o(\nu, \bar \eta )  + {\phi}_e(\nu, \bar \eta).
\]
The integration over the mixed regions  ($|z| < a$,
$|w| > b$, for example) vanishes by orthogonality.
The inner term reads
\begin{eqnarray}
\label{inner}
{\phi }_{o}(\nu,\bar \eta)
& = & \sum_{\ell=1}^{n-1} \frac{1}{(\nu \bar{\eta})^{\ell+1}} \sum_{m=n-\ell}^{n-1} c_m \int_{|z| < a} |z|^{2(m+\ell)}  d \mu_n(z) \\
                     &   &  + \sum_{\ell=n}^{\infty} \frac{1}{
                     (\nu \bar{\eta})^{\ell+1}} \sum_{m=0}^{n-1} c_m \int_{|z| < a} |z|^{2(m+\ell)}  d \mu_n(z) \nonumber   \\
                     &   &  + \sum_{\ell=0}^{n-1} \frac{1}{(\nu \bar{\eta})^{\ell+1}} \sum_{m=0}^{n- \ell - 1}
                                    c_m  c_{m+\ell}  \, \int_{|z| < a} |z|^{2(m+\ell)}  d \mu_n(z)  \, \int_{|z| > b} |z|^{2(m+\ell)}
                                    d \mu_n(z)   \nonumber  \\
                      & =  & { \phi}_{o}^{(1)}(\nu, \bar \eta    ) +     {\phi }_{o}^{(2)}(\nu, \bar \eta )
                                          +    {\phi}_{o}^{(3)}(\nu, \bar \eta).          \nonumber
\end{eqnarray}
The result of the outer integration has the same shape, up to the obvious inversions.  For example,
the analog of the first term on the right of (\ref{inner}) is  ${\phi}_e^{(1)}  = $
 $\sum_{\ell=1}^{n-1}  \sum_{m=n-\ell}^{n-1}  $
 $ (\nu  \bar{\eta})^{\ell} $
 $c_m \int_{|z| > b}  |z|^{2(m-\ell-1)} d \mu_n(z)$.
 We sketch the proof of the $L^1$ decay of $\phi_o$.
 The estimates for $\phi_e$ are much the same,  and both cases are similar to
 considerations immediately above.

Integrating by parts produces the bounds,
\begin{equation}
\label{gbdd0}
  \frac{n}{k+1} \, {a^{2k+2}} \, e^{-n a^2 }  \le    \int_{|z|< a}  |z|^{2k}  d \mu_n(z)  \le    \frac{1}{1 - \frac{n}{k}a^2}  \, \frac{n}{k+1}
    \, {a^{2k+2}} \, e^{-n a^2 },
\end{equation}
for $k > n$.  Recall here that $a < 1$.  If we further assume that $b > a$,
there is a constant $C$ such that
\begin{eqnarray*}
\label{gbdd1}
 \Bigl|     {\phi}_{o}^{(1)}(\nu, \bar \eta )   \Bigr|
& \le &   C \, e^{-n a^2}      \sum_{\ell = 0}^{n-1}   \sum_{m=0}^{\ell-1}
                   \frac{n^{n-(\ell - m)}}{(n- (\ell - m))!}  a^{2 ( n - (\ell - m))}  \, (a/b)^{\ell}    \\
&  \le &   C\,     e^{-n a^2}  \sum_{k = 0}^{n-1}  \frac{(na^2)^k}{k!}
                (a/ b)^{n-k},                  \nonumber
\end{eqnarray*}
(after  changing variables and the order of summation) and
\begin{equation*}
\label{gbdd2}
 \Bigl|      {\phi}_{o}^{(2)}(\nu, \bar \eta )     \Bigr|
 \le  C    \sum_{\ell = n}^{\infty}  (a/b)^{\ell}  \Bigl( e^{-na^2} \sum_{k=0}^{n-1}   \frac{(na^2)^k}{k!} \Bigr)
  \le  {C}    \frac{( a/b )^{n}}{  (1 - a/b)}.
\end{equation*}
The last two displays clearly tend to zero as $n \ra \infty$, providing
the advertised  $L^1$-decay over sets
supported  away from $|\nu| = |\eta|$.  To finish it is enough that
both
${\phi}_{o}^{(1)}(\nu, \bar \eta )$  and  ${\phi}_{o}^{(2)}(\nu, \bar \eta )$
remain bounded along that part of the  circle $|\nu| = |\eta|$ where $\nu /\eta = e^{i \theta}$
and say $|\theta| > \delta/2$.    This is a simple exercise taking advantage
of the oscillations introduced in summing powers of $e^{i \theta}$

As for  $ {\phi}_{o}^{(3)}(\nu, \bar \eta ) $,  the estimate (\ref{gbdd0})
does not have the same effect as $m + \ell \le n$
in each appearance of $\int  |z|^{2(m+\ell)} d \mu_n(z)$.
Instead, one works along the lines of (\ref{probbound}).
Again, first take $a < b$ and write,
\[
 \Bigl|      {\phi}_{o}^{(3)}(\nu, \bar \eta )     \Bigr| \le  \sum_{\ell = 1}^{n-1} (ab)^{-\ell} \sum_{m=0}^{n-\ell-1}
        \frac{c_{m+\ell}}{c_{m}} P ( S_{m+\ell}^n \le a^2 ) P ( S_{m+\ell}^n \ge b^2).
\]
Focussing on the sum restricted to  $\ell + m  \le na^2$, that object is bounded by
\begin{eqnarray*}
% \sum_{m+ \ell \le Na^2} (ab)^{-\ell} \frac{c_{m+\ell}}{c_{m}} P ( S_{m+\ell}^N \le a^2 ) P ( S_{m+\ell}^N \ge b^2)
   P(S_{na^2} \ge b^2)  \,  \sum_{\ell= 0}^{na^2}   (ab)^{-\ell}  \sum_{m=0}^{na^2 -\ell} \frac{c_{m+\ell}}{c_{m}}
  & \le &  2 e^{-n (b-a)^2/4} \,   \sum_{\ell= 0}^{na^2}   (ab)^{-\ell}
    \sum_{m=0}^{n a^2 -\ell} \frac{ (m^{\ell} +  \ell^2 m^{\ell-1} )}{n^{\ell}} \\
   & \le &  C n e^{-n (b-a)^2/4} \sum_{\ell = 0}^{ n a^2} \ell (a/b)^{\ell} \ra 0,
\end{eqnarray*}
having once more used the standard large deviations estimate for exponential random variables.
The sums over $m+\ell \ge n b^2$ and $ n a^2 \le m+\ell  \le n b^2$ are handled by the same
procedure.  When $a=b$  but $ \nu/ \eta $ is kept away from $1$, the key observation are that
$ P( S_{m+\ell}^n \le a^2 ) P ( S_{m+\ell}^n \ge a^2) $ has a good decay away from $| m + \ell - na^2| \le C \sqrt{n}$.
By using this in conjunction with the oscillations from $( \nu/\eta)^{-\ell}$, the boundedness
of $  {\phi}_{o}^{(3)}(\nu, \bar \eta )$ along $\{ |\nu| = |\eta| \} \cap \{ | \nu - \eta| > \delta \}$ will
follow.
\end{proof}

\begin{proof}[Proof of Lemma \ref{boundary}]   Start with
 the point-wise limit of $\phi_n(\nu, {\bar \eta})$ for which
we fix $\nu \neq \eta$ in the annulus $1 < |\cdot| < 1+\ep$. For the
interior estimate, it was be convenient to carry the integration in
$z$ and $w$ in the definition of $\phi_n$ over all of $\CC$, see
(\ref{phidef}).  However, by the estimate (\ref{ginbound}), one may
cut down the integral to $|\cdot| < 1 + \dl$  for any $\dl >0$ at
the expense of an exponentially small error. We therefore consider
the $n \ra \infty$ limit of
\begin{eqnarray*}
   {\widetilde \phi}_n(\nu, \bar{\eta})   & =  &
   \int_{|z| <  1 + \dl} \frac{1}{(z-\nu)(\overline{z - \eta})} K_n(z,\bar z) d \mu_n(z)  \\
      &  & -
\int_{ |z|, |w| < 1 + \dl } \frac{1}{(z-\nu)(\overline{w - \eta})}
                     | K_n(z,\bar w) |^2   \,  d \mu_n(z) d \mu_n(w),
\end{eqnarray*}
in which $\dl$ is chosen so that $1+\dl$ is less
than either $|\nu|$ or $|\eta|$.

Next, each appearance of $( \cdot -\eta)^{-1}$
and $(\cdot - \nu)^{-1}$ in ${\widetilde \phi}_n$ is expanded in series, and  the
integrals are performed  term-wise to find that,
\begin{eqnarray*}
   {\widetilde \phi}_n(\nu, {\bar{\eta}})  & =  & \sum_{\ell=0}^{\infty}  \int_{|z| < 1+\dl}
                                         \frac{|z|^{2 \ell}}{(\nu \bar{\eta})^{\ell+1}}  K_n(z,\bar{z})   d \mu_n(z) \\
                                   &      &  \  \ - \sum_{\ell  = 0}^{n-1}  \int_{|z| < 1+\dl} \int_{|w| < 1+\dl}
                                                          \frac{(z \bar{w}) ^{k}}{(\nu \bar{\eta})^{\ell+1}}  | K_n(z,\bar{w}) |^2
                                                                             d \mu_n(z)   d \mu_n(w) \\
                                                       &   =  &      \sum_{\ell =0}^{\infty}     \frac{1}{ (\nu \bar{\eta})^{\ell+1}}
                                                            \sum_{m=0}^{n-1}  c_m \int_{|z| < 1+\dl} |z|^{2(m+\ell)}  d \mu_n(z) \\
                                                       &        &  \  \ -     \sum_{\ell  = 0}^{n-1} \sum_{m=0}^{n-\ell -1}
                                                                               c_m  c_{m+\ell}  \Bigl(  \int_{|z| < 1+\dl}    |z|^{2(m +\ell)}  d  \mu_n(z) \Bigr)^2 .
\end{eqnarray*}
In the first equality, only the diagonal terms survive in the expansion an account of orthogonality.
Now note that, for any $\al > 0$,
\[
    \sum_{\al n \le  \ell  \le \infty}  \int_{|z| < 1+\dl}
                                             \frac{|z|^{2 \ell}}{|\nu{\eta}|^{\ell+1}}  K_n(z,\bar{z})  \mu_n(dz) \le n \sum_{\ell =
                                             \lfloor \al n \rfloor }^{\infty}  (1 - \theta_{\dl})^{\ell}
                                              = n \theta_{\dl}^{-1} (1 - \theta_{\dl})^{ \al n} ,
\]
in which  $\theta_{\dl} > 0 $ for $\dl > 0$, and
\[
     1 - c_k \int_{|z| < 1 +\delta} |z|^{2k} d  \mu_n(z)  =  O(e^{-n (\dl - \al)^2 } )
\]
for any $k \le (1+\al) n$ with $\al < \dl$.
It follows that,
\[
     {\widetilde \phi}_n(\nu, {\bar \eta})   =   \sum_{\ell=1}^{\lceil \al n \rceil}  \frac{1}{ (\nu \bar{\eta})^{\ell+1}}
                                                                   \sum_{m = n-\ell}^{n-1}  \frac{c_m}{c_{m+\ell}}    \,  +  o(1),
\]
with $\al$ obeying both the above constraints.   Finally, since
\[
    \sum_{m = n-\ell}^{n-1}  \frac{c_m}{c_{m+\ell}} =
         \sum_{m=0}^{\ell-1}  \Bigl(1 + \frac{m}{n} \Bigr) \Big(1 + \frac{m- 1}{n} \Bigr)
         \cdots \Bigl( 1 + \frac{ m - (\ell-1)}{n} \Bigr) =  \ell +  O \Bigl( \frac{\ell^3}{n^2}  \Bigr)
\]
we have that
\[
    \lim_{n \ra \infty}  \phi_n(\nu, {\bar \eta})  =  \sum_{\ell = 1}^{\infty}   \frac{\ell}{(\nu \bar{\eta})^{\ell+1}}.
\]
The limit  holds uniformly  on $\{   |\nu| \ge 1 + \delta \cup |\eta|  \ge  1 + \delta \}$ for any $\dl > 0$, and a dominated convergence
argument produces
\[
    \lim_{n \ra \infty}   \frac{1}{\pi^2} \int_{1 < |\nu|, |\eta| < 1+\ep }      {\bar \partial} f(\nu)  {\overline {{\bar \partial} f (\eta)}}  \phi_n(\nu, {\bar \eta})  \, d^2 \nu d^2 \eta
   =   \sum_{\ell = 1}^{\infty}   {\ell}  \,    \Bigl|   \frac{1}{\pi} \int _{1 < |\eta| <  1+\ep}   \frac{{\bar \partial}   f(\eta) }{\eta^{\ell+1} } \, d^2 \eta  \Bigr|^2.
\]
It remains to realize that since ${\bar \partial} [ f(\eta) \eta^{-k} ( \eta - \omega)  ]  = \eta^{-k} (\eta - \omega)  ({\bar \partial} f )(\eta)$ for
$\eta$ away from zero and any $\omega$, the dbar formula reads,
\[
   0 =  -  \frac{1}{\pi} \int _{1 < |\eta| <  1+\ep}   \frac{{\bar \partial}   f(\eta) }{\eta^{\ell+1} } \, d^2 \eta
     + \frac{1}{2 \pi } \int_{|\eta| = 1} \frac{ f(\eta)}{\eta^{\ell +1} }  \, d \eta  +  \frac{1}{2 \pi}   \int_{|\eta| = 1+\ep} \frac{ f(\eta)}{\eta^{\ell +1} }  \, d \eta.
\]
The third term vanishes by assumption, and the second term
is precisely $ \sqrt{-1}$ times the $\ell$-the Fourier coefficient
of the boundary data of $f$.
\end{proof}

\section{Universality for analytic functionals}
\label{sec:analytic}

%Here is an alternate proof of Theorem \ref{secondthm}, and more.  Once again it is profitable
%to take a slightly more general view, and consider a class of determinantal processes  in the complex
%plane of projection type with $n < \infty$ points.   In particular,
For each positive integer $n$, pick a rotation invariant measure,
$
      d \mu_n(z)  =   \tilde{\mu}_n( |z| ) d^2 z,
$
and set
$$
     K_n(z,\bar{w}) =  \sum_{\ell=0}^{n-1} c_{n,\ell} (z  {\bar w})^{\ell},  \
     \mbox{  where  }   \  c_{n,\ell}^{-1}= M(n,2\ell)= \int_{\CC}  |z|^{2 \ell}  d \mu_n(z);
$$
the measure $\mu_n$ itself need not depend explicitly on $n$.
Now impose the following moment condition:
for all integer $m$ there
 is a $\rho>0$ such that
\begin{equation}
\label{momentass} M(n,2n+2m)/M(n,2n)\to \rho^{2m},\qquad
\end{equation}
as $n \ra \infty$.  Of course,
by a scaling, it may be assumed that $\rho = 1$.

For any determinantal point process $(K, \mu)$ in $\CC$ with
rotation invariant $\mu$, it is the case that the moduli $\{
|z_1|, |z_2|, \dots, |z_n| \}$ as a set have the same distribution
as $\{R_1,\ldots, R_k\}$ where $R_i$ are independent with $P( R_k
\in dr ) = 2 \pi c_k r^{2k +1}$ $ \tilde \mu(r) dr$ for each $k$
(see \cite{BenKrPerVir05} or \cite{Kostlan92}). Therefore,
(\ref{momentass}) has the interpretation that the stochastically
largest random variable $R_n$ satisfies $\ev_n R_n^{2m} \to \rho$,
for every $m\in \Z$. Equivalently, for every fixed $\ell\ge 0$ we
have $\ev_n R_{n-\ell}^2\to \rho$.

Being the main example at hand, the reader will be happy to check
that Ginibre satisfies everything asked for: the relevant
computation is $\lim_{n \ra \infty}  \frac{\Gamma( n + m )}{n^m
\Gamma( n) } = 1$.  A second example of interest is the {\em
truncated Bergman ensemble}.  Here one begins with $\mu_n = $ the
uniform measure on $\UU$, producing the kernel
\[
   K_n(z, \bar {w}) = \sum_{\ell = 0}^{n-1} (\ell+1) z^{\ell} {\bar w}^{\ell}
\]
of orthonormal polynomials on the disk.  This model arises
naturally in the following way.  Consider the random polynomial
$z^n +  \sum_{k=0}^{n-1} a_k z^k$ with independent coefficients
drawn uniformly from the large disk $R\, \UU$. If we condition the
roots to lie entirely within the unit disk, the $R\to\infty$ limit
of the resulting point process is the truncated Bergman ensemble.
This observation may be gleaned from Hammersely \cite{Ham56}, but
see also \cite{PerVir05} for why the adjective ``truncated" is
used.

Consider now a linear statistic in any such ensemble which is polynomial in $z$ alone,
or of the form
$
\sum_{i=1}^n p_m(z_i) =
\sum_{k=1}^m a_k  \Bigl( \sum_{i=1}^n z_i^k \Bigr).
$
We have the following CLT.

\begin{theorem}
\label{momentthm} Let $z_1,\dots, z_n$ be drawn from a
determinantal process $(K_n , \mu_n)$ as above for which
(\ref{momentass}) holds with $\rho=1$. Denote $p_j(z^{\oplus
n})=z_1^j+\cdots +z_n^j$. Take $a = (a_1, \dots, a_k)$ and $b =
(b_1, \dots, b_k)$ with $a_j, b_j \in \{0,1,\dots \}$ and bounded
independently of $n$. Then,
\begin{equation}
\label{momenteq}
  \lim_{n \ra \infty}  \ev_n \left[  \prod_{j=1}^k  p_j(\zn)^{a_j} \overline{ p_j(\zn)}^{b_j} \right]
    =  \ev \left[ \prod_{j=1}^k  \Bigl(  \sqrt{j} Z_j \Bigr)^{a_j}   \overline{   \Bigl(  \sqrt{j} Z_j \Bigr)^{b_j} }  \right],
\end{equation}
for ${ Z}_1, \dots , {Z}_k$ independent standard complex Gaussian
random variables.
\end{theorem}

Having identified the limiting moments, this implies Theorem 3.  On a case by case basis,
Theorem \ref{momentthm} and a variance estimate yields a more complete picture.

\begin{corollary}
\label{analyticor}
Let $f(z)$ be analytic in a neighborhood of  $|z| \le  1$ and the
points $z_1, \dots, z_n$ be drawn either from the Ginibre or
truncated Bergman ensemble.  Then,
as $n \ra \infty$,
$
   \sum_{\ell = 1}^n f(z_{\ell})  - n f(0)
$
converges in distribution to a mean-zero complex Normal with variance
 $\frac{1}{\pi}  \int_{\UU} | f^{\prime} (z) |^2 \, d^2 z.$
\end{corollary}

Qualitatively, the Ginibre and truncated Bergman ensemble have marked
differences.  While the Ginibre eigenvalues fill the disk as $n \ra \infty$, the
truncated Bergman points concentrate near $|z| = 1$
({\em i.e.},  $\frac{1}{n} K_n(z, \bar z)$ tends weakly to $\delta_{|z| = 1}$).
Once again, the analytic
CLT only ``sees the boundary''.
%Moreover, Theorem \ref{momentthm} shows the
%limiting noise to be  identical for a large class of rotation invariant determinantal
%processes, at least if one restricts to polynomial statistics.

The proof  is largely inspired by the ideas of Diaconis-Evans
\cite{DiacEv01} where the analogous moment formula is established
for the eigenvalues of the Haar distributed unitary group. Those
points yield yet another example in the above class: there $d
\mu_n(z) = \delta_{|z| = 1} $ and $\rho = 1$.  Also in that case,
\cite{DiacEv01} shows the equality (\ref{momenteq}) to hold for
finite $n$ as soon as $n \ge ( \sum_{j=1}^k j a_j ) \vee (
\sum_{j=1}^k j b_j )$. The basic observation is that the integrand
on the left hand side of (\ref{momenteq}) is comprised of symmetric
polynomials in the points $z_1, \dots, z_n$; one would like to to
expand  this object in a convenient basis.

Let $A_n$ be the vector space of symmetric polynomials in the
variables $z^{\oplus n}=z_1,...,z_k$ of degree at most $n$. Let
$\Lambda_n$ denote the set of partitions of integers at most
$n$.  Given a partition $\lambda=(\ld_1 \ge \ld_2 \ge \cdots
\ge ...) \in \Lambda_n $, bring in the
corresponding {\em Schur function},
\be
\label{defschur}
   s_{\ld}( \zn) =  \frac{ \det \Bigl(  z_{k}^{\ld_{\ell} +n - {\ell}} \Bigr)_{k, \ell = 1}^n }
                                                          { \det \Bigl( z_k^{n - \ell}  \Bigr)_{k, \ell  = 1}^n },
\ee
and recall the following well known facts.

\begin{theorem}\label{schur} (\cite{Mac79}, Chapter 1)
The Schur functions $s_{\ld}( z^{\oplus n}), \lambda\in
\Lambda_n$ form a basis for $A_n$. Consider the inner product
$\langle\cdot, \cdot \rangle$ that makes them an orthonormal
basis. For $f_n(\zn)=\prod_{j=1}^k p_k^{a_k}(\zn)$ and
 $g_n( \zn)  =\prod_{j=1}^k p_k^{b_k}(\zn)\in A_n$, which are  products
of simple power sum functions, it holds that
 \be
\label{inner2} \langle f_n,g_n \rangle  = \delta_{ab}  (
\prod_{j=1}^n j^{a_j}a_j! ).
 \ee
The inner products are compatible as $n$ varies in the sense
that if $m>n$, and $\lambda\in \Lambda_n$ we have $\langle
f_n, s_\lambda(\zn) \rangle=\langle f_m,s_\lambda(\zm) \rangle$.
\end{theorem}

\begin{lemma}
\label{schurcomp} For Schur functions in the points of $(K,
\mu_n)$ as above, \be \label{eschur}
   \ev_n \Bigl[  s_{\ld} (\zn) \, \overline{s_{\pi} (\zn) } \Bigr]
   = \delta_{\lambda\pi} \prod_{\ell=0}^{n-1} \frac{M(n,2(n+\ld_{\ell}
   -\ell))}{M(n,2(n-\ell))}.
\ee
\end{lemma}

\begin{proof}
The denominator in \re{defschur} is a Vandermonde determinant, and
thus the interaction term (see \re{dens}) in the integral of
question is cancelled:
\[
  \ev_n \Bigl[  s_{\ld}( \zn) \, \overline{s_{\pi}(\zn) } \Bigr]
 =    \frac{1}{{\mathcal Z}_n} \,
  \int_{\CC} \cdots \int_{\CC} \,
    \det \Bigl( z_k^{\ld_{\ell} + n - \ell } \Bigr)  \det \Bigl(  z_k^{\pi_{\ell}  + n - \ell} \Bigr)
      d \mu_n(z_1) \dots d \mu_n(z_n). \nonumber
\]
The normalizer here is  ${\mathcal Z}_n = n! \prod_{\ell = 1}^n
M(n,2\ell)$. Now, expanding each determinant on the right hand
side, the generic term we get is a constant multiple of
\[
       \prod_{k=1}^n  \int_{\CC}   z^{\ld_{k} - k}   {\bar z}^{\pi_{\sigma_k }- \sigma_k}   |z|^{2n}      \, d \mu_n(z)
\]
with a permutation $\sigma$.   Since $\ld$ and $\pi$ are
monotone, this will vanish by orthogonality of $z$ and $\bar z$
when $\ld \neq \pi$. When $\ld=\pi$, there is a non-zero
contribution
from the diagonal, $\sigma = \it{id}$. To
conclude, we compute
\[
   \ev_n  \Bigl[  s_{\ld}(\zn) \, \overline{s_{\ld}(\zn) } \Bigr]
   =   \frac{n!}{{\cal Z}_n} \,  \prod_{k=0}^{n-1}  \int_{\CC} |z|^{2(\ld_{\ell} + n- \ell)} d
   \mu_n(z),
\]
but this is exactly (\ref{eschur}).
\end{proof}

\begin{proof}[Proof of Theorem \ref{momentthm}]
Let $f_N=\prod_{j=1}^k  p_j(\zn)^{a_j}$ and $g_n=
\prod_{j=1}^k p_j(\zn)^{b_j}$, and let $n_0$ be at
least of the degrees of $f$ and $g$ (this is independent of $n$).
Now expand $f_n,g_n$ with respect to the basis given by the
Schur functions, and compute the expectation term by term. This
expansion is finite, and the coefficients $\langle f_n,
s_\lambda(\zn) \rangle$, $\langle g_n,
s_\lambda(\zn) \rangle$ do not depend on $n$ for $n \ge
n_0$ by Theorem \ref{schur}.

Since \re{eschur} has a bounded number of factors that are not 1,
our assumption and Lemma \ref{schurcomp} implies that  for
$\lambda, \nu$ fixed we have $\ev_n [ s_\lambda  \bar{s_\nu}  ]
\to \delta_{\lambda \nu}$. Thus we have $\ev_n [f_n \bar g_n] \to
\langle f_{n_0},g_{n_0}\rangle$. The claim now follows by
\re{inner2} and the moment formula for complex Gaussians.
\end{proof}

\begin{proof}[Proof of Corollary \ref{analyticor}]  For analytic $f(z)$,
$f(0) = \frac{1}{\pi} \int_{\UU} f(z) d^2 z =  \frac{1}{2 \pi}   \int_{ \partial \UU} f(z) dz$, which explains
the centralizer.  While for the Ginibre ensemble, we could simply quote the result of Section \ref{sec:conc},
the analyticity allows for a simpler approach amenable to more general ensembles subject  perhaps to
additional conditions in the spirit of (\ref{momentass}).  Preferring to be concrete (and brief),
we restrict ourselves to the truncated Bergman case.

Rather than using the dbar representation, the variance in the analytic case may be computed by the more
familiar Cauchy integral formula:  with the reference measure now uniform on the unit disk,
\[
    \frac{1}{4 \pi^2} \int_{C_{\delta}}  \int_{C_{\delta}}  f(\nu) \overline{f(\eta)}   \Bigl[   \int_{\UU}
        \frac{K_n(z,\bar z)}{(z -\nu)( \bar{z}  - \bar{\eta} )} d^2 z
      -  \int_{\UU} \int_{\UU}   \frac{|K_n(z,\bar w)|^2}{(z -\nu)( \bar{z}  - \bar{\eta} )}   d^2 z d^2 w  \Bigr]  d \nu d  \bar{\eta}.
\]
Here $C_{\delta}$ is a circle of radius $1 + \delta > 1 $ (about the
origin)  within the region of analyticity of $f$.
What to do  next is plain:
\[
  \int_{\UU}
        \frac{K_n(z,\bar z)}{(z -\nu)( \bar{z}  - \bar{\eta} )} d^2 z
= \sum_{\ell = 0}^{\infty} \frac{1}{(\nu \bar \eta)^{\ell + 1} }
          \sum_{m=0}^{n-1} \frac{m+1}{m+ \ell +1},
\]
and
\[
  \int_{\UU}  \int_{\UU}   \frac{|K_n(z,\bar w)|^2}{(z -\nu)( \bar{w}  - \bar{\eta} )}   d^2 z d^2 w
    =  \sum_{\ell = 0}^{n-1} \frac{1}{(\nu \bar \eta)^{\ell + 1} }
          \sum_{m=0}^{n - \ell -1} \frac{m+1}{m+ \ell +1},
\]
after expanding both the kernel  and the functions $(z - \cdot)^{-1}$, $(w - \cdot)^{-1}$ and
integrating term-wise.
With $|\nu \eta | > 1$ there is enough control to pass the limit inside the summations
and conclude that the variance tends
to
\begin{eqnarray*}
    \frac{1}{4 \pi^2}  \int_{C_{\delta}}  \int_{C_{\delta}}  f(\nu) \overline{f(\eta)}  \Bigl( \sum_{\ell=1}^{\infty}  \frac{\ell}{(\nu \bar \eta)^{\ell + 1} }
  \Bigr)   d \nu d  \bar{\eta} & = & \sum_{\ell=1}^{\infty} \ell  \, \Bigl|  \frac{1}{2 \pi}  \int_{\partial U}  \frac{f(\nu)}{\nu^{\ell +1}} d \nu  \Bigr|^2,
\end{eqnarray*}
as $n \ra \infty$. Here we have used the analyticity to pull the
integral back to $\partial \UU$ after the limit was performed. This
last expression, and much of  the method getting there, is now
recognized from the Ginibre boundary case (Lemma \ref{boundary}). In
particular, that it equals the advertised $ \frac{1}{\pi}\int_{\UU}
| f^{\prime}(z)|^2 d^2 z$ has already been explained.
\end{proof}

\bigskip

\noindent {\bf{Acknowledgements.}}  The work of B.R. was supported
in part by NSF grant DMS-0505680,  and that of B.V. by a Sloan
Foundation fellowship, by the Canada Research Chair program, and
by NSERC and Connaught research grants. Both authors thank the
hospitality of MSRI (spring  2005) where  this project was
initiated.

%\noindent{\bf Acknowledgements.} We thank
%\bibliography{}

\sc \bigskip \noindent B\'alint Vir\'ag, Departments of
Mathematics and Statistics, University of Toronto, ON, M5S 2E4,
Canada.
\\{\tt balint@math.toronto.edu}, \ {\tt
www.math.toronto.edu/\~{}balint}

\sc \bigskip \noindent B. Rider, Department
of Mathematics, University of Colorado at Boulder,
Boulder, CO 80309. \\{\tt brider@euclid.colorado.edu},
\ {\tt math.colorado.edu/\~{}brider}
\end{document}